\documentclass[12pt]{article}
\usepackage{latexsym}
\usepackage{amsmath}
\usepackage{amsfonts}
\usepackage{amssymb}
\usepackage{amsthm}
\usepackage{amscd}
\usepackage{color}
\input xy
\xyoption{all} \hfuzz=1.5pt \tolerance=500

\newcommand{\ot}{\otimes}

\newcommand{\va}{\varphi}
\newcommand{\Om}{\Omega}

\newcommand{\id}{{\bf 1}}
\newcommand{\co}{{\B C}}
\newcommand{\s}{\sigma}
\newcommand{\lm}{\lambda}

\newcommand{\prtpet}{\mathbin{\widehat{\otimes}}}

\newcommand{\bd}{\begin{document}}
\newcommand{\ed}{\end{document}}
\newcommand{\Htp}{\mathbin{\stackrel{\cdot}{\bigotimes}}}
\newcommand{\Hts}{\mathbin{\stackrel{\cdot}{\bigoplus}}}
\newcommand{\nl}{\nu\in\Lambda}
\newcommand{\xl}{X_\nu;\nl}
\newcommand{\x}{X_\nu}
\newcommand{\la}{\langle}
\newcommand{\ra}{\rangle}
\newcommand{\ch}{\cal H}
\newcommand{\su}{\subseteq}
\newcommand{\B}{\Bbb}
\newcommand{\e}{\varepsilon}
\newcommand{\ov}{\overline}
\newcommand{\vk}{\varkappa}
\newcommand{\vt}{\vartheta}
\newcommand{\al}{\alpha}
\newcommand{\de}{\delta}
\newcommand{\lo}{\Longleftrightarrow}
\newcommand{\cK}{{\cal F}}
\newcommand{\sq}{{\square}}
\newcommand{\bb}{{\cal B}}
\newcommand{\dd}{{\cal D}}
\newcommand{\cR}{{\cal R}}
\newcommand{\f}{{\cal K}}
\newcommand{\cb}{{\cal CB}}
\newcommand{\wrr}{\widetilde{\cal R}}
\newcommand{\br}{^\bullet{\cal R}}
\newcommand{\en}{E_\nu}
\newcommand{\el}{{\cal L}}
\newcommand{\lr}{\Longrightarrow}
\newcommand{\bp}{\bigcirc}
\newcommand{\Long}{\Longleftarrow}
\newcommand{\q}{\quad}
\newcommand{\qq}{\qquad}
\newcommand{\cd}{\cdot}
\newcommand{\fur}{f: E\to{\B R}}
\newcommand{\furo}{f_0: E_0\to{\B R}}
\newcommand{\fuc}{f: E\to{\B C}}
\newcommand{\ii}{\infty}
\newcommand{\di}{\diamondsuit}
\newcommand{\bgd}{{\bigtriangledown}}
\newcommand{\bu}{{\bigtriangleup}}
\newcommand{\bc}{{completely bounded}}
\newcommand{\cc}{{completely contractive}}
\newcommand{\qs}{{quantum space}}
\newcommand{\isc}{{isometric}}
\newcommand{\ism}{{isomorphism}}
\newcommand{\qss}{{quantum spaces}}
\newcommand{\bco}{{completely bounded operator}}
\newcommand{\bcos}{{completely bounded operators}}
\newcommand{\res}{{respectively}}
\newcommand{\tp}{{tensor product}}
\newcommand{\eq}{{equivalent}}
\newcommand{\qtp}{{quantum tensor product}}
\newcommand{\mma}{\mathrel{\mathop{\otimes}\limits_{A}}}
\newcommand{\mmb}{\mathrel{\mathop{\otimes}\limits_{\bb}}}
\newcommand{\mmp}{\mathrel{\mathop{\otimes}\limits_{p}}}
\newcommand{\mmh}{\mathrel{\mathop{\otimes}\limits_{h}}}
\newcommand{\mmop}{\mathrel{\mathop{\otimes}\limits_{4}}}
\newcommand{\mmf}{\mathrel{\mathop{\otimes}\limits_{op}}}
\newcommand{\mmi}{\mathrel{\mathop{\otimes}\limits_{i}}}
\newcommand{\mmm}{\mathrel{\mathop{\otimes}\limits_{\bb-\bb}}}
\newcommand{\mmdd}{\mathrel{\mathop{\otimes}\limits_{\cdot}}}

\newcommand{\mmo}{\mathrel{\mathop{\cdot}\limits_{1}}}
\newcommand{\mmA}{\mathrel{\mathop{\otimes}\limits_{A-A}}}
\newcommand{\mmd}{\mathrel{\mathop{\cdot}\limits_{2}}}
\newcommand{\msp}{\mathrel{\mathop{\otimes}\limits_{sp}}}
\newcommand{\mms}{\stackrel{h}{\otimes}}
\newcommand{\mmt}{\stackrel{4}{\otimes}}
\newcommand{\mmx}{\stackrel{i}{\otimes}}
\newcommand{\mmy}{\stackrel{p}{\otimes}}
\newcommand{\mmz}{\stackrel{sp}{\otimes}}
\newcommand{\gd}{\ddagger}
\newcommand{\od}{\odot}
\newcommand{\mt}{\mapsto}
\newcommand{\mmc}{\mathrel{\mathop{\otimes}\limits_{\sim}}}
\newcommand{\mme}{\stackrel{\sim}{\otimes}}

\begin{document}

{\centerline{{\bf Metric freedom and projectivity for }}

{\centerline{{\bf classical and quantum normed modules }}

   \vspace{1cm}

\centerline{A.~Ya.~Helemskii}

\centerline{Faculty of Mechanics and Mathematics}

\centerline{Moscow State University}

\centerline{Moscow 119991 Russia}

\vspace{1cm}

{\bf Abstract}

\medskip
In functional analysis there are several reasonable approaches to the notion of a projective
module. 
%Some of them take into account the norm topology of the modules in question, and some (like
%a new notion of the metric projectivity) the exact value of the norm.
%At the beginning 
We show that a certain general-categorical framework 
%suggest a general-categorical framework that 
contains, as particular cases, all known versions. In 
this scheme, the notion of a free object comes to the forefront, and in the best of categories, 
called freedom-loving, all projective objects are exactly retracts of free objects. We concentrate 
on the so-called metric version of projectivity and characterize metrically free `classical', as 
well as quantum ( = operator) normed modules. Hitherto known the so-called extreme projectivity 
turns out to be, speaking informally, a kind of `asymptotically metric projectivity'. 

Besides, we answer the following concrete question: what can be said about metrically
projective modules in the simplest case of normed spaces ? We prove that metrically projective
normed spaces are exactly $l_1^0(M)$, the subspaces of $l_1(M)$, where $M$ is a set, consisting of finitely
supported functions. Thus in this case the projectivity coincides with the freedom.

\section{ Introduction}

 There are several essentially different ways to carry over to functional analysis the notion of a
projective module, one of most important in homological algebra. Four of them seem to be most 
important, two taking into account the norm topology of the modules in question, and  two depending 
on the exact value of the norm. (The comparatively recent interest to the latter two was stimulated 
by some questions of operator space theory; cf.~\cite{he4}\cite{he5}\cite{wit}). Here are the 
relevant definitions. 

\medskip
Let ${\cal K}$ be a category, so far arbitrary, $P,Y,X$ its objects, $\tau:Y\to X$, $\va:P\to X$ its morphisms. We recall that a morphism $\psi:P\to Y$ is called a \emph{lifting of $\va$ across $\tau$}, if it makes the diagram
$$
\xymatrix@R-10pt@C+15pt{
& Y \ar[d]^{\tau}\\
P \ar[ur]^{\psi} \ar[r]^{\va} & X }\eqno(1.1)
$$
\noindent commutative.

\medskip
Throughout the paper, $A$ is an arbitrary contractive normed algebra, supposed, for simplicity, to be 
unital.
Saying `normed module' or just `module', we mean a contractive left normed unital $A$-module. The 
identity map on a set $M$ will be denoted by $\id_M$.

We begin with the oldest and most known variety, introduced (in the context of Banach modules) in~\cite{he0}; see also, 
e.g.,~\cite{he1}~\cite{he2}\cite{run}.

\medskip
{\bf Definition 1.1}. A module $P$ is called {\it relatively projective} (or just {\it projective},
as in the great majority of literature), if for every bounded module morphism $\tau:Y\to X$ such
that there exists a bounded operator $\rho:X\to Y$ with $\tau\rho=\id_X$, every bounded morphism 
$\va:P\to X$ has a bounded lifting across $\tau$.

\medskip
The following, more demanding condition cropped up in connection with the study of amenable Banach algebras
(cf.~\cite{he1}). However, in the context of Banach spaces ( = Banach $\co$-modules) it was known long ago under the name of the lifting property (see,  e.g.,~\cite{pie} and the references therein).

%  These are exactly spaces,
%topologically isomorphic to $l_1(M)$ for some index set $M$. (see ???).

\medskip
{\bf Definition 1.2}. A module $P$ is called {\it topologically projective}, if, for every bounded 
surjective module morphism $\tau:Y\to X$, every bounded morphism $\va:P\to X$ has a bounded lifting 
across $\tau$. 

\medskip
If $E$ is a normed space,  we denote by $\bp_E$ its closed unit ball, and by $\circledcirc_E$ its open unit ball . Recall that an operator $\tau:E\to F$ between normed spaces is called {\it coisometric} (respectively, {\it strictly coisometric}), if it maps $\circledcirc_E$ onto
$\circledcirc_F$ (respectively, $\bp_E$ onto $\bp_F$). Coisometries are often called quotient maps
and strict coisometries exact quotient maps. Recall that the Hahn/Banach Theorem admits the
following equivalent formulation:

-- \emph{an operator $\va : E\to F$ between normed spaces is an isometry if, and only if its
adjoint $\va^*:F^*\to E^*$ is a coisometry, and if, and only if this adjoint is a strict
coisometry}. (This will be crucial to much of what will follow).

\medskip
{\bf Definition 1.3}. A module $P$ is called {\it extremely projective}, if, for every coisometric
module morphism $\tau:Y\to X$ and every $\e>0$, every bounded morphism $\va:P\to X$ has a lifting $\psi$ across $\tau$ such that $\|\psi\|<\|\va\|+\e$.

\medskip
Extremely projective Banach spaces are known in the Banach space geometry as `spaces with the
metric lifting property'~\cite[p. 133]{pie}. Initially these spaces were considered by
%in the well known paper of
Grothendieck~\cite{gro}. (What he has done, will be recalled and used below). Much
later, extremely projective normed
 and Banach modules were
formally defined and studied in~\cite{he7}.
%, where some related results were obtained.
(See~\cite{he6} about the so-called extremely flat and
extremely injective modules). Why we say `extreme', see Remark 2.8 below.

\medskip
The following notion is new, and it is one of main topics of this paper. As we shall see, in some questions (in particular, from the general-categorical overview) it seems to behave better, than the extreme projectivity.

\medskip
{\bf Definition 1.4}. A module $P$ is called {\it metrically projective}, if, for every strictly
coisometric  module morphism $\tau:Y\to X$, every bounded morphism $\va:P\to X$ has a lifting
$\psi$ across $\tau$ such that $\|\psi\|=\|\va\|$.

\medskip
Every kind of projectivity has a non-completed, `normed' as well as the completed, `Banach' 
version. (The relevant behaviour of modules in both cases is often very similar, but sometimes 
essentially different; cf.~\cite[Thm.  ]{he7}). Besides, each of these versions has a `classical' 
prototype and its `quantum' (= operator space) counterpart. In particular, the quantum extreme 
projectivity (under the name of just projectivity) was introduced and studied by Blecher~\cite{ble} 
in 1992. 

Note a useful link between `non-completed' and `completed' versions.

\medskip
{\bf Proposition 1.5.} {\it Let $P$ be a normed module over a Banach algebra $A$, $\ov{P}$ its completion. Suppose that $P$ possesses a projectivity property, formulated in one of definitions
1.1 -- 1.4, in its `non-completed' version.   Then $\ov{P}$ possesses the same property in its `completed' version.}

\smallskip
PROOF. We restrict ourselves with the case of the extreme projectivity. Suppose we are given a coisometric morphism $\tau:Y\to X$ between Banach modules, a bounded morphism
$\va:\ov{P}\to X$ and $\e>0$. Consider the restriction $\va_0$ of $\va$ to $P$. It has a
lifting $\psi_0:P\to Y$ of $\va_0$ across $\tau$ with $\|\psi_0\|<\|\va_0\|+\e$. Since $Y$ is complete, $\psi_0$ has the continuous extension $\psi:\ov{P}\to Y$. The rest is clear. $\Box$

\bigskip
The present paper pursues several aims. First, we show that a certain general-categorical framework 
%suggest a general-categorical framework that 
contains, as particular cases, all 
%known versions. In suggest a general-categorical scheme that 
%contains, as particular cases, all 
mentioned versions of projectivity as well as their injective 
counterparts. The basic notion is the so-called rigged category, generalizing the relative abelian 
category of MacLane~\cite[Ch.IX.5]{mc1} (Here we mostly concentrate on categories that are not even 
additive). Within this framework, the notion of a free object comes to the forefront. In fact, in 
this paper  we study projective objects by means of free objects. 

 In the majority of examples such a scheme is sufficient to describe
 projective objects as retracts of free objects. Similarly,
injective objects are
%and their description as
retracts of the so-called cofree objects.

A notable exception is the case of the extreme projectivity (see above). The `asymptotic' nature of such a
notion, with the indispensable $\e$ in its definition, requires a kind of elaboration of our
scheme. This is achieved by supplying the given rig by an additional, the so-called asymptotic
structure. In this framework, the notion of the asymptotically projective object naturally
appears, and such an object can be characterized as the so-called asymptotic retract of a free object (Proposition 6.12). This general scheme includes, as particular cases, extremely projective classical, as well as quantum modules.

\bigskip
Our second aim is to give an explicit description of free (and cofree) objects in our main examples 
of rigged categories of modules. We describe metrically free and cofree `classical' modules. 
Further, we find a suitable rig of the category of quantum modules and characterize relevant free 
objects, demonstrating their abundance. In particular, in the case of the simplest base algebra 
$\co$ the free quantum module, that is just the free quantum space, with the one-point base turn 
out to be of the form 
$$
{\cal N}_1\oplus_1{\cal N}_2\oplus_1\cdots\oplus_1{\cal N}_n\oplus_1\cdots,
$$
where ${\cal N}_n$ is the space of nuclear, or trace class operators on $\co^n$, equipped with the trace
quantum norm, and $\oplus_1$ is the symbol of the quantum $l_1$-sum of quantum spaces. (See details in Section 5).

We emphasize that the part of our paper, dealing with quantum spaces and modules, is written under the strong influence of the already mentioned paper~\cite{ble}.
 Among other ideas and results, Blecher was first to realize
the crucial role of spaces ${\cal N}_n$ in respective lifting problems, and he has described, in terms of
these spaces, all extremely projective quantum spaces. (See Corollary 6.13(ii) and relevant comments).

\bigskip
The third (last but not the least) aim of this paper is to give a full description of metrically projective normed, not necessarily completed spaces. (We answer the similar question about Banach spaces as well, but this is easier).
%, that is modules over the simplest base algebra $\co$. Theorem 3.?
We show that the former are exactly
the spaces $l_1^0(M)$, the normed subspaces of $l_1(M)$, consisting of finitely supported
functions; here $M$ is an index set.
This means, in particular, that in the context of normed (as well as of Banach) spaces the metric projectivity coincides with the metric freedom. We do not know whether such a coincidence is
true for the class of {\it extremely} projective normed spaces, which is, generally speaking, larger.
(Note that in the context of Banach spaces the answer is `yes', thanks to Grotendieck; see details in Section 3).

\bigskip
The order of our presentation does not coincide with the given enumeration of our aims; in this way we tried to make the paper more readable. In Section 2 we introduce projective and free objects as derived notions of that of a rigged category, and then characterize free objects in the context of `classical' normed and Banach modules (Theorem 2.18).

In Section 3 we leave, for a time, `abstract nonsense' and deal with fairly `concrete'
question of the characterization of metrically projective normed spaces (Theorem 3.5).

The subjects of Section 4 are injectivity and cofreedom. We introduced these notions, again in the framework of rigged categories, and characterize cofree `classical' modules with the help of some general-categorical observations (Theorem 4.9 and Proposition 4.5).

In Section 5 we suggest the rigged category, appropriate to deal with
metrically projective quantum spaces, and we describe the relevant free objects (Theorem 5.9).

Finally, in Section 6 we show that the notion of the extremely projective module in classical and quantum
contexts, as well as results of the characterization of these modules,
fit to the general framework of an asymptotic category (Corollary 6.13).

\section{ Projectivity and freedom in rigged categories. `Classical' examples}

\bigskip
\q\;{\bf Definition 2.1.} Let ${\cal  K}$ be an arbitrary category. A {\it rig} of ${\cal  K}$ is a faithful
%(see, e.g.,~\cite[p. 15]{mc2})
covariant functor $\square:{\cal K}\to{\cal L}$, where ${\cal L}$ is another category. A pair, consisting of a category and its rig, is called {\it rigged category}.

\medskip
Fix, for a time, a rigged category, say $({\cal K},\square:{\cal K}\to{\cal L})$. We call a
morphism $\tau$ in ${\cal K}$ {\it admissible epimorphism}, if $\square(\tau)$ is a retraction in
${\cal L}$. (Clearly, such a $\tau$ is indeed an epimorphism.)

\medskip
{\bf Definition 2.2}. An object $P$ in ${\cal K}$ called {\it $\square$-projective} (or, if there
is no danger of misunderstanding, (just) \emph{projective}), if, for every admissible epimorphism
\\ $\tau:Y\to X$ and every morphism $\va:P\to X$, there exists a lifting of $\va$ across $\tau$.

In a concise form, $P$ is projective, if the standard covariant morphism functor ${\bf h}_{{\cal
K}}(P,?):{\cal K}\to{\bf Sets}$ takes admissible epimorphisms to surjective maps.

\medskip
Let us distinguish an immediate (and actually well known)

\medskip
{\bf Proposition 2.3}. {\it (i) A retract of a projective object is itself projective

(ii) If $\s:X\to P$ is an admissible epimorphism in ${\cal K}$, and $P$ is projective, then $\s$ is
a retraction.} $\Box$

\bigskip
Turn, for a time, to some examples. In what follows, $A$ is an
% given normed
algebra (cf. Introduction),
%Further,
${\bf A-mod}$ is the category of normed
$A$-modules and their bounded morphisms, and ${\bf A-mod_1}$ is
the category of the same modules and their \emph{contractive} morphisms. By ${\bf Nor}$ and ${\bf
Nor}_1$ we denote, as usual, the category of normed spaces ( = ${\co}$-modules) and all bounded,
respectively, all contractive operators. The corresponding categories of Banach ( = complete)
modules over Banach algebras and of Banach spaces we denote by  $\ov{{\bf A-mod}}$, $\ov{{\bf
A-mod_1}}$, ${\bf Ban}$ and ${\bf Ban}_1$. As usual, ${\bf Set}$ denotes the category of sets.

\bigskip
{\bf Example 2.4.} Set ${\cal K}:=\ov{{\bf A-mod}}$, ${\cal L}:={\bf Ban}$ and take as $\sq$ the
respective forgetful (about the outer multiplication) functor. Then, of course, we obtain the old
and well known definitions of an admissible epimorphism of Banach $A$-modules and of a (relatively)
projective Banach $A$-module. (see, e.g.,~\cite{he1}~\cite{he2} or~\cite{run}). The obvious
non-completed version of the mentioned rig, with ${\cal K}:={{\bf A-mod}}$ and ${\cal L}:={\bf
Nor}$ leads to the notion
%s of an admissible epimorphism of normed $A$-modules and
of a (relatively)
projective normed $A$-module, that is to the kind of projectivity, mentioned in Definition 1.1.

\bigskip
{\bf Example 2.5.} Now ${\cal K}$ is again $\ov{{\bf A-mod}}$, but ${\cal L}$ is ${\bf Set}$; so
$\sq$ `forgets about everything'. Then we obtain all surjective (and hence open) morphisms of
Banach $A$-modules in the capacity of admissible epimorphisms, and the obvious completed version of
topologically projective modules from Definition 1.2 in  the capacity of $\sq$-projective objects.
These modules were explicitly mentioned in~\cite{he7}; their flat and injective counterparts
appeared under another name many years ago in~\cite{he1}.

In the `non-completed' case we replace in the latter rig
$\ov{{\bf A-mod}}$ by ${{\bf A-mod}}$ and get the kind of projectivity, introduced in
Definition 1.2. However, in this case admissible morphisms are not be bound to be open, and, as a
result, we shall obtain too few projective objects. For example, in the
case $A:=\co$ projective normed spaces are only those that are finite-dimensional.

\medskip
 {\bf Remark 2.6.} This nuisance can be removed. E.~Gusarov~\cite{gus} has shown that, if
we shall consider, as a rig, the forgetful functor from ${{\bf A-mod}}$ into the category ${\bf
Bor}$ of the so-called bornological sets, then the resulting admissible morphisms will again be
open and the family of projective objects considerably increases.
%In particular, normed spaces that are projective relatively to this rig, can be described
%THANKS TO GROENBAEK,
% as topologically isomorphic to $l_1^0(M)$ for some index set $M$.

\bigskip
We turn to the main rig of this section.

\medskip
{\bf Example 2.7.} Consider the functor
%Set ${\cal K}:={\bf A-mod_1}$, ${\cal L}:={\bf Set}$, and take as $\sq$ the functor
$\bigcirc:{\bf A-mod_1}\to{\bf Set}$, taking a module $X$  to its {\it closed} unit ball, and a 
contractive morphism  to its birestriction to the closed unit balls of the respective modules. 
Obviously, in this case \emph{admissible morphisms are exactly strictly coisometric morphisms} and, 
accordingly, projective modules are metrically projective modules from Definition 1.4. The 
indicated rig has an obvious version for Banach modules. 

\medskip
{\bf Remark 2.8.} It is easy to observe that metrically projective objects can be also defined in
the following way: a normed $A$-module $P$ is metrically projective, if the morphism functor ${\bf
h}_A(P,?):{\bf A-mod}\to{\bf Nor}$ preserve the property of a morphism to be strictly coisometric.
Similarly, $P$ is extremely projective, if this functor
 preserve the property of a morphism to be (just) coisometric, or, equivalently, to be an extreme
epimorphism in the general-categorical sense (cf. ~\cite[1.7]{clm} or~\cite[p. 100]{he3}. (Hence the word `extreme').

\medskip
{\bf Remark 2.9.} If we shall consider the obvious analogue of the rigged category in Example 2.7
by taking the symbol $\circledcirc$ (see Introduction) instead of $\bp$, our admissible morphisms
will turn out to be (just) coisometric morphisms. However, we shall not get, in the capacity of
projective modules, extremely projective modules from Definition 1.3; moreover, as it is actually
well known, we shall have zero module as the only $\circledcirc$-projective module. We do not know
a rig that would provide extremely projective modules, and our conjecture is that such a rig does
not exist. Nevertheless, as it was mentioned in Introduction, there is a device, allowing one to
put the extreme projectivity on the formal basis. It will be discussed in Section 6.

\medskip
We restricted ourselves with the mentioned examples of rigs. Note only that there is a lot of
others, and some of them can be rather curious. One of them is the rig $\ov{{\bf A-mod}}\to{\bf
Ban}^{\bf o}$, where ` $^{\bf o}$ ' denotes the dual category; it takes a module to the dual of its underlying
normed space. This rig is useful because, compared with that in Example 2.4, it provides
much larger stock of admissible epimorphisms.

 \bigskip
Let us return to the general context and fix, for a time, a rigged category \\
$({\cal K},\square:{\cal K}\to{\cal L})$. The following concept is well known under different 
names.
% , despite one can acroat least in pure algebra. 

\medskip
{\bf Definition 2.10} 
%(cf.~\cite{tsal})~\cite[Ch.2.7]{faith} or~\cite[Ch.0.7]{he3}). 
Let $M$ be an object in ${\cal L}$. An object ${\bf 
Fr}(M)$ in ${\cal K}$ is called {\it free (or, to be precise, $\square$-free) object with the base 
$M$}, if,  for every $X\in{\cal K}$, there exists a bijection 
$$
{\cal I}_X: {\bf h}_{{\cal L}}(M,\square X)\to{\bf h}_{{\cal K}}({\bf Fr}(M),X),\eqno(2.1)
$$
natural in the second argument. A rigged category is called {\it freedom-loving}, if every object
in ${\cal L}$ is a base of a free object in ${\cal K}$.

\medskip
Note that in a freedom-loving rigged category the assignment $M\mt{\bf Fr}(M)$ can be extended to
morphisms in the appropriate `coherent' manner (cf.~\cite[p. 81]{mc2}). In this way we get the
so-called freedom functor ${\bf Fr}:{\cal L}\to{\cal K}$,
%. Besides, it is easy to show (cf. {\it idem}) that a freedom functor
which is exactly a left adjoint functor to $\sq$;
see, e.g.,\cite[Ch.IV.1]{mc2}. We recall that, for arbitrary functors $\Psi:{\cal K}\to{\cal L}$
and $\Phi:{\cal L}\to{\cal K}$, $\Phi$ is called left adjoint to $\Psi$ (or, equivalently, $\Psi$
is called right adjoint to $\Phi$), if, for every $X\in{\cal K}$ and $M\in{\cal L}$, there exists a
bijection
$$
{\cal I}_{M,X}: {\bf h}_{{\cal L}}(M,\Psi(X))\to{\bf h}_{{\cal K}}(\Phi(M),X),
$$
 natural in {\it both} arguments.

The following observations show the practical use of the freedom. They are actually well known and 
can be extricated, as particular cases or easy corollaries, from some general facts, contained 
in~\cite[Chs. III,IV]{mc2}. 

\medskip
{\bf Proposition 2.11}. {\it (i) If, for an object $X$ in ${\cal K}$, the object $\sq(X)$ in ${\cal
L}$ is the base of a free object, then there exists an admissible epimorphism of this free object
onto $X$. Thus in a freedom-loving rigged category every object in ${\cal K}$ is a range of an
admissible epimorphism with a free domain.

(ii) All free objects in ${\cal K}$ are projective.

(iii)If our rigged category is freedom-loving, then an object in ${\cal K}$ is is projective if,
and only if it is a retract of a free object.}

\medskip
As we shall see on concrete examples in this and especially in 5th section, the search of
projective and free objects can be considerably facilitated, if our categories admit coproducts.
Let us recall, in a concise form, what it is.

A \emph{coproduct} of a family $X_\nu;\nl$ of objects in ${\cal K}$ is a pair $(X,\{i_\nu;\nl\})$,
where $X$ is an object, and $i_\nu:X_\nu\to X$ are morphisms. By definition, this pair has the following property: for
every object $Y$ the map between ${\bf h}_{\cal K}(X,Y)$ and the cartesian product
$\textsf{X}\{{\bf h}_{\cal K}(X_\nu,Y);\nl\}$, taking $\psi$ to the family $\{\psi i_\nu\}$, is a
bijection. The mentioned $X$, denoted in the detailed form by $\coprod\{X_\nu;\nl\}$, will be
referred as the {\it coproduct object}, $i_\nu$ as the {\it coproduct injections}, and the indicated property
as the \emph{universal property of a coproduct}. See, e.g.,~\cite[Ch.2]{faith} or~\cite[p.
59]{he3}. The morphism $\psi$ will be called the {\it coproduct of the morphisms $\psi_\nu$.}

The definition easily implies that
if we have two coproducts of the same family of objects, then there is a categorical isomorphism
between the respective coproduct objects, compatible in an obvious way with the respective
coproduct injections. Therefore we shall speak about `the' coproduct of a given family.

We say that ${\cal K}$ \emph{admits coproducts}, if every family of its objects has the
coproduct.

Now take our given rigged category.

\medskip
{\bf Proposition 2.12.} \emph{Suppose that $(P,\{i_\nu;\nl\})$ is the coproduct of the family
$P_\nu;\nl$ of projective objects in ${\cal K}$. Then $P$ is also projective}.

\smallskip
PROOF. Take $\tau$ and $\va$ as in Definition 2.2. Consider the respective liftings
$\psi_\nu:P_\nu\to Y$ of morphisms $\va_\nu:=\va i_\nu$ and their
coproduct $\psi$. We see that $\tau\psi$, as well as $\va$, is the coproduct of morphisms
$\va_\nu;\nl$. Hence $\tau\psi=\va$. $\Box$
%these morphisms coincide. $\Box$
%, by the `uniqueness' part of the universal property, $\tau\psi=\va$. $\Box$
%\medskip
%Turn to the freedom.

\medskip
{\bf Proposition 2.13.} \emph{Suppose that for an index set $\Lambda$ and every $\nl$ we are given
a free object $F_\nu$ with a base $M_\nu$ Suppose, further, that there exist the coproducts
$F:=\coprod\{F_\nu;\nl\}$ in ${\cal K}$ and $M:=\coprod\{M_\nu;\nl\}$ in ${\cal L}$. Then $F$ is a
free object with the base object $M$.}

\smallskip

PROOF. Using the definition of a coproduct and that of a free object, consider the chain of
bijections
$$
{\bf h}_{{\cal L}}(M,\sq(X))=\textsf{X}\{{\bf h}_{\cal L}(M_\nu,\sq(X);\nl\}=
\textsf{X}\{{\bf h}_{\cal K}(F_\nu,X);\nl\}={\bf h}_{{\cal K}}(F,X),
$$
obviously natural in $X$. The rest is clear. $\Box$

%Using the definition of a coproduct, then that of a free object, and then again that of a
%coproduct, we obtain the bijection between ${\bf h}_{{\cal L}}(M,\sq(X))$ and the cartesian product
%of sets ${\bf h}_{{\cal L}}(M_\nu,\sq(X))$ , then between the mentioned cartesian product and the
%cartesian product of all ${\bf h}_{{\cal K}}(F_\nu,X)$ and finally between the latter cartesian
%product and ${\bf h}_{{\cal K}}(F,X)$. As the composition, we obtain a bijection between ${\bf
%h}_{{\cal L}}(M,\sq(X))$ and ${\bf h}_{{\cal K}}(F,X)$, obviously natural in $X$. $\Box$

\bigskip
Come back to  our concrete examples of rigged categories. Do they love freedom?

Throughout the paper, we use the notation `$\ot_p$' and `$\prtpet$' as the symbols of non-completed 
and completed projective tensor product, respectively. The spaces of the form $A\ot_p E$ or 
$A\prtpet E$ are considered as normed or, according to the sense, Banach $A$-modules with respect 
to the outer multiplication, well defined by $a\cd(b\ot x):=ab\ot x; a,b\in A, x\in E$ (cf., 
e.g.,~\cite[Ch. III.1]{he1}). 

\medskip
{\bf Example 2.14.} The rigged category from Example 2.4 is freedom-loving, and this is well known.
For every Banach space $E$, the free Banach left $A$-module with a base $E$ is $A\prtpet E$
(see~\cite{he1} or~\cite{he2}). The same is true for the non-completed version of the rigged
category in question, only free modules have the form $A\ot_p E$.

\medskip
{\bf Example 2.15.} At the same time, the rigged category in Example 2.5 is not freedom-loving. One 
can easily observe that a Banach $A$-module is free if, and only if it is topologically isomorphic 
to the module $A\prtpet\co^n$ for some $n=1,2,...$, and its base set consists of $n$ points. As a 
corollary, only finite sets can be base sets of free Banach modules, and the latter are finitely 
generated. 
%\smallskip
%PROOF. We shall only show that an infinite set, say $M$, can not be the base of a free module.
%Suppose that, on the contrary, that there is a certain ${\bf Fr}(M)$ with the base $M$ and a
%universal arrow $\al$. Take $Y:=A\prtpet l_1(M)$ and consider the map $\va^0:M\to\sq(Y)$, taking
%$t\in M$ to the vector $\lm_te_+\ot\de_t$ where $\lm_t>0$ is chosen in such a way, that
%$\inf\{\|\al(t)\|/\lm_t; t\in M\}=0$. Consequently, if an $A$-module morphism $\va:{\bf Fr}(M)\to
%Y$ makes the diagram (1.2) commutative, then we have $\va(\al(t)=\va^0(t)$ and hence
%$\|\va(\al(t)\|=\lm_t$. Therefore $\va$ can not be bounded, a contradiction. $\Box$

\medskip
{\bf Remark 2.16.} If we replace, in the capacity of ${\cal L}$, ${\bf Set}$ by ${\bf Bor}$
 (cf. Remark 2.6), we obtain much larger stock of free objects and their bases (E.~Gusarov~\cite{gus}). In
particular, Gusarov has shown that in such a context a module is free if, and only if it is
topologically isomorphic to $A\prtpet l_1(M)$ for some set $M$.

\medskip
Needless to say, the category ${\bf Set}$ admits coproducts (the so-called disjoint unions of given
sets). Therefore Proposition 2.16 has an immediate

\medskip
{\bf Corollary 2.17.} \emph{Suppose that ${\cal K}$ admits coproducts, ${\cal L}$ is just ${\bf Set}$,
and there exists a free module, say $F$, with one-point base. Then our rigged category is
freedom-loving. Moreover, if $M$ is an arbitrary set, then the coproduct object of the family of
copies of $F$, indexed by points of $M$, is the free object with the base $M$.}

\medskip
In particular, this corollary shows that the poverty of the stock of free objects in Example 2.15 
is connected with the following known limitation of the category $\ov{{\bf A-mod}}$: the family of 
non-zero objects in this category has a coproduct if, and only if this family is finite (cf., 
e.g.,~\cite[Ch.2.5]{he3}). 

On the contrary, we recall that categories ${\bf A-mod_1}$ and $\ov{{\bf A-mod_1}}$ \emph{admit 
coproducts}. Indeed, it is easy to see that the coproduct object of the family $X_\nu;\nl$ in ${\bf 
A-mod_1}$, respectively $\ov{{\bf A-mod_1}}$, is the non-completed, respectively, completed 
$l_1$-sum of the given modules, together with the natural embedding of direct summands. Recall also 
that the non-completed $l_1$-sum of a family of copies of the module $A$, indexed by the points of 
a set $M$, is not other thing than $A\ot_p l_1^0(M)$, whereas the completed $l_1$-sum of such a 
family is $A\prtpet l_1(M)$. 
%(This is because of Grothendieck's identification of Banach spaces $E\prtpet l_1(M)$ and $l_1(M,E)$,
%where $E$ is a Banach space, and its obvious non-completed counterpart).
In particular, the space $l_1^0(M)$, respectively, $l_1(M)$ is the coproduct object in ${\bf
Nor_1}$, respectively, in ${\bf Ban_1}$ of the family of copies of $\co$, indexed by points of $M$.

Armed with these facts, we return to our main rig $\bigcirc:{\bf A-mod_1}\to{\bf Set}$. In what
follows, $\bigcirc$-free objects of ${\bf A-mod_1}$ will be referred as \emph{metrically free
modules}.

\medskip
{\bf Theorem 2.18}. {\it The rigged category $({\bf A-mod_1}, \bigcirc)$ is freedom-loving.
Moreover, the metrically free normed $A$-module with a set $M$ as its base is $A\ot_p l_1^0(M)$.}

\smallskip
PROOF. By what was said above,
%virtue of Corollary 2.17 and the above-mentioned description of coproducts in ${\bf A-mod_1}$,
we must only show that, for a one-point set, say $\{t\}$, $A$ is a
metrically free module with this base. But this is indeed the case: for every object
$X$ in ${\bf A-mod_1}$, the map between ${\bf h}_{\bf Set}(\{t\},\bigcirc(X))$ and ${\bf h}_{{\bf
A-mod_1}}(A,X)$, taking a map $\va$ to the (contractive) morphism $\psi:a\mt a\cd\va(t)$, is
obviously a bijection, natural in $X$. $\Box$

%Indeed, define $\al:\{t\}\to\bp(A)$ by taking $t$ to $e_+$. Then we immediately see that, for a
%given left normed $A$-module $Y$ and a map $\va^0:\{t\}\to\bp(Y)$, the morphism $\va:A\to Y:a\mt
%a\cd\va^0(t)$ is a unique morphism, making the respective specialization of the diagram (2.1) (with
%$\al$ as $\al_M$ etc.) commutative. The rest is clear. $\Box$

\medskip
This theorem has an obvious Banach (= completed) version. Namely, metrically free Banach 
$A$-modules with the base $M$ are those of the form $A\prtpet l_1(M)$. 

\medskip
{\bf Corollary 2.19.} \emph{Every normed, respectively, Banach $A$-module is  the image of the
module $A\ot_p l_1^0(M)$, respectively, $A\prtpet l_1(M)$ for some set $M$ under a strictly
coisometric morphism. Besides, such a module is metrically projective if, and only if for some set
$M$ it is a retract (or, equivalently, module direct summand with the natural projection of norm 1)
of a module  $A\ot_p l_1^0(M)$, respectively, $A\prtpet l_1(M)$}.

\medskip
Of course, the mentioned facts have simple direct proofs, but we wanted to show that they, as well
as much else, are but particular manifestations of a general scheme we discuss.

\bigskip
{\bf Remark 2.20.} On the contrary, it is easy to show that the rigged category $({\bf A-mod_1},
\circledcirc:{\bf A-mod_1}\to{\bf Set})$ (cf. Remark 2.9) is not freedom-loving, and, moreover,
there is no $\circledcirc$-free $A$-modules save 0.}

%\smallskip
%PROOF. Suppose that, on the contrary, there exists a $\circledcirc$-free module, say $F$, with the
%base $M$ and universal arrow $\al_M:M\to\circledcirc(F)$. Fix an arbitrary point $t\in M$. Of
%course, there exists a vector $x\in\circledcirc(F)$ with $\|\al_M(t)\|<\|x\|$. Consider the map
%$\va^0:M\to\circledcirc(F)$, taking all $M$ to $x$. Since $F$ is free, we must have a morphism
%$\va:F\to F$ in ${\bf A-mod_1}$, taking $\al_M(t)$ to $x$. But morphisms in our present category
%can not increase norms of vectors. We came to a contradiction. $\Box$

\bigskip
{\bf Remark 2.21.} Many years ago, Semadeni~\cite{sem0} proposed an essentially different 
definition of free, as well of a projective objects in a large class of categories. His approach 
fits very well to a lot of important examples, but, according to what is said in ({\it idem}, pp. 
5,27), not to categories of modules over general algebras. 
%, more complicated than $\co$.

\section{ Metrically projective normed spaces are $l_1^0$}

%We interrupt, for a time, the discussion about metric and other kinds of projectivity from the
%overview of the `abstract nonsense'  and
Now we turn to a quite concrete question of the geometry of
normed spaces. What is the structure of metrically projective spaces ?

Recall that Grothendick actually answered a similar question for {\it extremely} projective {\it Banach} spaces:
they happened to be $l_1(M)$, where $M$ is a set. It is natural to conjecture that extremely
projective normed spaces are $l_1^0(M)$, but we do not know whether it is true. However, we can
describe {\it metrically} projective spaces, non-completed and (what is easier) completed alike.
This is what will be done in this section. We begin with some preparation.

We call a Banach space, say $E$, {\it metrically flat} if, for every isometry $i:F\to G$ of Banach spaces
the operator $\id_E\prtpet i:E\prtpet F\to E\prtpet G$ is also isometry.
%Similarly, in terms of non-completed tensor product ` $\mmp$ ', we define a metrically flat normed space.

\medskip
{\bf Proposition 3.1.} {\it Every metrically or extremely projective Banach space is  metrically flat.}

\smallskip
PROOF. Let $E$ be a given space, $i:F\to G$ as above. We recall
%By the equivalent formulation on The Hahn-Banach Theorem
(see Introduction) that the property of the operator $\id_E\prtpet i$ to be
isometric is equivalent to the property of its adjoint $(\id_E\prtpet i)^*$ to be strictly
coisometric as well as to be (just) coisometric. By the adjoint associativity (now in its simplest
form; see, e.g.,~\cite[p. 180]{he3} or~\cite[(6.1.4)]{he4}), this adjoint is isometrically
equivalent to the operator $\bb(E,i^*):\bb(E,G^*)\to\bb(E,F^*):\va\mt i^*\va$. But, since $i$ is
isometric, $i^*$ is strictly coisometric and hence coisometric. Therefore, by the assumption on
$E$, $\bb(E,i^*)$ is coisometric. $\Box$

\medskip
In our subsequent argument we turn to Grotendieck's theorem, mentioned above. Note that in
literature, speaking about this theorem, they usually cite~\cite{gro} (see, e.g.,~\cite[p.
182]{heb}). Basically, it is correct. At the same time, despite the paper~\cite{gro} contains all
needed ingredients for the proof, the theorem itself is not explicitly formulated.
 By `ingredients' we mean the following two statements, formulated (needless to say, in equivalent terms)
 and completely proved:

\medskip
(i) A Banach space is metrically flat if and only if it is isometrically isomorphic to some
$L_1(\Omega,\mu)$, where $(\Omega,\mu)$ is a measure space~\cite[Prop. 2]{gro}

(ii) A closed subspace of $l_1(N)$, where $N$ is a set, topologically isomorphic to some
$L_1(\Omega,\mu)$, is isometrically isomorphic to $l_1(M)$ for some set $M$~\cite[Prop. 2]{gro}.

\medskip
{\bf Proposition 3.2.} {\it Every metrically projective Banach space is isometrically isomorphic to
$l_1(M)$ for some set $M$.}  

\smallskip
PROOF. By the assumption,
%Being metrically projective,
our space is a retract in ${\bf Ban_1}$ of some free Banach
space, that is of $l_1(N)$ for some set $N$ (Corollary 2.19). Hence it coincides, up to an isometric
isomorphism, with a closed subspace of $l_1(N)$. At the same time, combining Proposition 3.1 with
(i) above, we see that our space is isometrically isomorphic to some $L_1(\Om,\mu)$. Then  the
assertion (ii) above works. $\Box$

\medskip
{\bf Remark 3.3.} We see that the same argument, using Proposition 3.1, completes the proof of the
Grothendieck Theorem itself as it was formulated above. (The only difference is that we use,
instead of Corollary 2.19, Proposition 6.13(i) below). Apparently, the traditional way to prove this
proposition is to use a non-trivial criterion of the metric flatness, namely Proposition 1
in~\cite{gro}, and then to verify the relevant condition. For this aim one takes a certain family
of operators (indexed by $\e$ from Definition 1.3), and then, applying the Banach-Alaoglu Theorem,
proceeds to the cluster point of this family with respect to a suitable weak$^*$ topology. See,
e.g.,~\cite[27.4.2]{sem}. The way we suggest above seems to be shorter.

\medskip
However, the non-completed version of Proposition 3.2 needs some additional work. We came to the
main result of this section.

\medskip
{\bf Theorem 3.5.} {\it Every metrically projective normed space is isometrically isomorphic to
$l_1^0(M)$ for some set $M$.}

\medskip
PROOF. We see from Propositions 1.5 and 3.2 that our given space, let it be $E$, is, up to an
isometric isomorphism, a dense subspace of some $l_1(M)$. Therefore our task is to show
that our space contains all vectors of $l_1^0(M)$ and nothing more.

Combining Proposition 2.11(iii) and Theorem 2.18, we see that $E$ is a retract in ${\bf Nor}_1$ of
$l_1^0(N)$ for some set $N$. Fix a retraction $\s:l_1^0(N)\to E$ and a coretraction
 $\rho:E\to l_1^0(N)$ with $\s\rho=\id_E$. Consider also the respective extentions by continuity
$\ov{\s}:l_1(N)\to l_1(M)$ and $\ov{\rho}:l_1(M)\to l_1(N)$. For a set $L$, we denote by $bas(L)$
the set of characteristic functions of one-point subsets of $L$; this set is, of course, a linear basis
in $l_1^0$.

Take an arbitrary $e\in bas(M)$. Evidently, $\ov{\rho}(e)$ has a unique expansion $\sum_k\lm_ke_k'$
with no more than countable set of summands, such that vectors $e_k'$ are multiples of pairwise
different vectors of $bas(N)$, $\|e_k'\|=1$ and $\lm_k>0$ for all $k$.
%Note that $e_k'$ lies in $l_1^0(N)$, that is in the domain of $\s$.

\medskip
{\bf Lemma}. {\it We have $\s(e_k')=e$ for all $k$. }

\smallskip
Fix some $e_k'$. Since $\|e\|=1$ and $\ov{\rho}$, together with ${\rho}$, is an isometry, we have
$\sum_k\lm_k=1$. Of course, we can suppose that the expansion of $\ov{\rho}(e)$ has at least two
summands. It follows that $\ov{\rho}(e)$ is a convex combination of $e_k'$ and $z:=(\sum_{l:l\ne
k}\lm_l)^{-1}\sum_{l:l\ne k}\lm_le_l$.
Consequently, $e=\ov{\s}\,\ov{\rho}(e)$  is a convex combination of $\ov{\s}(e_k')$ and
$\ov{\s}(z)$. Since $\ov{\s}$, as well as $\s$, is contractive, both of these vectors belong to
$\bp_{l_1(M)}$. But $e$ is, of course, an extreme point of $\bp_{l_1(M)}$. Therefore
$\ov{\s}(e_k')=e$. It remains to recall that $e_k'$ lies in the domain of $\rho$.

\medskip
{\bf The end of the proof.}

\medskip
The lemma shows that
%$bas(M)$ lies in $E$, and consequently
$E$ contains all $l_1^0(M)$.
Now suppose that $E$ contains vectors, not belonging to $l_1^0(M)$. Then there exists $x\in S_E$ of
the form $\sum_{k=1}^\ii\lm_ke_k$, where, for all $k\in{\B N}, e_k$ are pairwise different vectors
of $bas(M)$, and $\lm_k\ne0$.

Take, for every $k$, the
%, e_k$ in the capacity of $e$ above and consider the
expansion of $\rho(e_k)$ of the form $\sum_l\lm_{kl}e_{kl}'$, where all $e_{kl}'$ are multiples of
different vectors of $bas(N)$, $\|e_{kl}'\|=1$ and $\lm_{kl}>0$ for all indexes.
%respective standard expansion of $\rho(e_k)$; let it be $\sum_l\lm_{kl}e_{kl}'$.
By Lemma, we have $\s(e_{kl}')=e_k$ for all $e_{kl}'\in bas(N)$, participating in the latter sum.
It obviously follows that for different $k$ the supports of $\rho(e_k)$ (as of functions on $N$) do
not intersect. Consequently, the support of $\rho(x)=\sum_{k=1}^\ii\lm_k\rho(e_k)$ is the disjoint
union of the supports of $\rho(e_k); k\in{\B N}$. Therefore it is an infinite set. But, on the
other hand, $\rho(x)$ lives in $l_1^0(N)$ and hence its support is finite, a contradiction.  $\Box$

\medskip
{\bf Remark 3.6.} As it was mentioned, we do not know whether {\it extremely} projective normed 
spaces are isometrically isomorphic to $l_1^0(M)$; our method, using extreme points of unit balls, 
does not work. At the same time, the similar question about {\it topological} projectivity is now 
answered. Groenbaek~\cite{groe} has shown that every topologically projective normed space is {\it 
topologically} isomorphic to $l_1^0(M)$. This is the non-completed version of an earlier result of 
K\"othe~\cite{heb} who has proved that every topologically projective {\it Banach} space is 
topologically isomorphic to $l_1(M)$. 

\section{ Injectivity and cofreedom}

We return to the general scheme.
 What if one wishes to introduce, together with projective objects or instead of them, injective
objects in a given category ${\cal K}$ ? To provide for them a formal basis, consider a rigged
category $({\cal K},\boxdot:{\cal K}\to{\cal M})$, where ${\cal M}$) is an auxiliary category
(generally speaking, different from ${\cal L}$ in Section 2).

We call a morphism (necessarily monomorphism) $\iota$ in ${\cal K}$ {\it admissible}, if
$\boxdot(\iota)$ is a coretraction in ${\cal M}$. Then we call an object $J\in{\cal K}$ {\it
injective} (to be precise, \\ $\boxdot$-injective), if the standard contravariant morphism functor
${\bf h}_{{\cal K}}(?,J):{\cal K}\to{\bf Set}$ takes admissible monomorphisms  to surjective maps.

The same can be said in the concise form with the help of the so-called rigged category, dual to
$({\cal K},\boxdot)$. This is the pair $({\cal K}^{\bf o},\boxdot^{\bf o}:{\cal K}^{\bf o}\to{\cal M}^{\bf o})$,
where $^{\bf o}$ is a symbol of the dual category as well as of the relevant `copy` of the functor
$\boxdot$. Namely, $J$ is injective with respect to the initial rigged category exactly if it is
projective with respect to the dual rigged category.

\medskip
{\bf Example 4.1.} Consider the rigged category $(\ov{{\bf A-mod}},\sq:\ov{{\bf A-mod}}\to{\bf
Ban})$ from Example 2.4. It is easy to see that we obtain the traditional notions of an admissible
monomorphism between Banach $A$-modules and a (relatively) injective Banach $A$-module (see
again~\cite{he1}~\cite{he2}~\cite{run}).

\medskip
{\bf Example 4.2.} Consider the rigged category $\circleddash:{\bf A-mod}\to{\bf Set}^{\bf o}$, 
taking a module $X$ to the underlying set of its dual module $X^*$, and a morphism in ${\bf A-mod}$ 
to its adjoint, considered as just map. Using the Hahn/Banach Theorem, one can easily see that 
$\circleddash$-admissible are  topologically injective morphisms (cf., e.g.,~\cite[Ch. 2.5]{he3}). 
As to $\circleddash$-injective Banach modules, we call them topologically injective. The latter, 
under the name of strictly injective modules, play a certain role in Banach 
homology~\cite[Ch.VII.1]{he1} (although lesser role than relatively injective modules). 
%Obviously, topologically injective Banach spaces (= ``{\B 0}-modules'') are exactly those
%possessing the so-called extension property (cf., e.g.,~\cite[p. 133]{pit}).

%A similar procedure, applied to the forgetful functor from Example 2, provides topologically
%injective morphisms in the capacity of admissible monomorphisms and the so-called topologically
%injective Banach modules. The latter, under the name of strictly injective modules, where mentioned
%in~\cite[44]{hel4}. Obviously, topologically injective Banach space (= ``0-modules'') are exactly
%those possessing the so-called extension property (cf., e.g.,~\cite[44]{pit}).

\medskip
{\bf Example 4.3.} Now, as a counterpart to the rigged category in Example 2.7, consider the rigged
category $({\bf A-mod_1},\circledast:{\bf A-mod_1}\to{\bf Set}^{\bf o})$, where $\circledast$ is a
covariant functor, taking a module $X$ to the closed unit ball $\circledast_X:=\bigcirc_{X^*}$ of
its dual module $X^*$, and a morphism in ${\bf A-mod_1}$ (which is, as we remember, contractive) to the
respective restriction of its adjoint to unit balls.

It follows from the equivalent formulation of the Hahn/Banach Theorem (see Introduction) that
$\circledast$-admissible monomorphisms are exactly the isometric morphisms, and
$\circledast$-injective objects are exactly those with the `Hahn/Banach property' (`metric
injective property', as they say in Banach space geometry). So, it is justified to call these
modules {\it metrically injective}.

\medskip
{\bf Example 4.4.} If we replace, in the definition of the latter rig, the ball $\bigcirc_{X^*}$ by
$\circledcirc_{X^*}$, nothing will change. The
%above-mentioned equivalent formulation of the
Hahn/Banach Theorem provides the same isometries as admissible monomorphisms and hence the same
injective objects.

%In order to define what would be called ``extremely injective modules'' (cf. Example 4), it seems
%natural to replace in the definition of the functor $\circledust$ the word ``closed'' by ``open''.
%However, nothing will change. It easily follows from the Hahn/Banach Theorem (more precisely, from
%its equivalent form, mentioned before), that the respective admissible monomorphisms are the same
%as in the previous example.

\medskip
Turning from the pro-- to injectivity, we inevitably come to the so-called cofreedom. We say, for
brevity, that the object, say ${\bf Cfr}(M)$, in ${\cal K}$ is cofree with respect to the given
rigged category $({\cal K},\boxdot:{\cal K}\to{\cal M})$ , and $M\in{\cal M}$ is its cobase, if
${\bf Cfr}(M)$ is free with respect to the dual rigged category $({\cal K}^{\bf o},\boxdot^{\bf o})$, and $M$
is its base. In a
%n obvious
way, parallel to Definition 2.10, we define a {\it cofreedom-loving} category.
We
%immediately
see that Proposition 2.11 has its appropriate counterpart; in particular, \emph{an
object in a cofreedom-loving category is injective if, and only if it is a retract of a cofree object.}

The following general-categorical observation provides a unified method to describe cofree objects
in a lot of concrete cases.

%It is not hard to describe these modules, using the special property of the rig in question.
%However, we shall do this as a particular case of a certain general observation.

\medskip
{\bf Proposition 4.5.} Let $({\cal K}_1,\boxdot:{\cal K}_1\to{\cal L}_1)$, $({\cal K}_2,\sq:{\cal
K}_2\to{\cal L}_2)$ be rigged categories, and
%supplied by covariant functors
$\Psi:{\cal L}_1\to{\cal L}_2$, $\Upsilon:{\cal K}_1\to{\cal K}_2$, covariant functors, making the
diagram
$$
\xymatrix@C+20pt{{\cal K}_1 \ar[r]^{\boxdot}\ar[d]_{{\Psi}}
& {\cal L}_1 \ar[d]^{\Upsilon} \\
{\cal K}_2 \ar[r]^{\sq} & {\cal L}_2 }\eqno(4.1)
$$
\noindent commutative. Suppose that $\Psi$ and $\Upsilon$ have left adjoints $\Phi$ and $\Delta$,
respectively, and $F$ is a free object in ${\cal K}_2$ with the base $M$. Then $\Phi(F)$ is a free
object in ${\cal K}_1$ with the base $\Delta(M)$.

\smallskip
PROOF. Take an arbitrary object $Y$ in ${\cal K}_2$ and consider the chain of bijections
$$
{\bf h}_{{\cal L}_1}(\Delta(M),\boxdot(Y))\to{\bf h}_{{\cal L}_2}(M,\Upsilon\boxdot(Y)\to
{\bf h}_{{\cal L}_2}(M,\sq\Psi(Y)\to
$$
$$
{\bf h}_{{\cal K}_2}(F,\Psi(Y))\to{\bf h}_{{\cal K}_1}(\Phi(F),Y),
$$
provided by the assumption on $\Delta$ and $\Phi$, the definition of a free object and the diagram
(4.1). The resulting bijection between  ${\bf h}_{{\cal L}_1}(\Delta(M),\boxdot(Y)$ and ${\bf
h}_{{\cal K}_1}(\Phi(F),Y)$ is obviously natural in $Y$. The rest is clear. $\Box$

\medskip
To apply this proposition to our principal examples, consider, together with the category ${\bf
A-mod}$, its `right-module twin' ${\bf mod-A}$. Every normed space $E$ gives rise to the standard
contravariant functors $\bb(?_l,E):{\bf A-mod}\to{\bf mod-A}$ and $\bb(?_r,E):{\bf mod-A}\to{\bf
A-mod}$ (cf., e.g.,~\cite[Ch.III.1]{he1}). Here $\bb(\cd,\cd)$ means a relevant space of all bounded operators,
equipped with the operator norm. The first functor takes $X$ to $\bb(X,E)$ with
the right outer multiplication, defined by $(T\cd a)(x):=T(a\cd x);a\in A, x\in X, T\in\bb(X,E)$
and $\va:X\to Y$ to $\bb(\va,E):\bb(Y,E)\to\bb(X,E):\psi\mt\psi\va$. The second functor
is defined in a similar pattern.

In the same way, using the notation ${\bf mod-A_1}$ for the appropriate categories of right
modules, we obtain contravariant functors from ${\bf A-mod_1}$ to ${\bf mod-A_1}$, and from ${\bf
mod-A_1}$ to ${\bf A-mod_1}$. For them we retain the notation $\bb(?_l,E)$ and
$\bb(?_r,E)$.

\medskip
{\bf Proposition 4.6.} \emph{The functor $\bb(?_r,E)$, being considered as a covariant functor from
${\bf mod-A}$ to ${\bf A-mod}^{\bf o}$, is a left adjoint to the functor $\bb(?_l,E)$, the latter being
considered as a covariant functor from ${\bf A-mod}^{\bf o}$ to ${\bf mod-A}$. The same assertion
holds, if we replace ${\bf mod-A}$ by ${\bf mod-A_1}$ and ${\bf A-mod}^{\bf o}$ by ${\bf A-mod_1}^{\bf o}$.}
%, up to obvious modifications in its formulation, holds for the functors $\bb(?_r,E): ${\bf
%mod-A_1}\to{\bf A-mod_1}^{op}$ and $\bb(?_l,E):{\bf A-mod_1}^{op}\to{\bf mod-A_1}$.}}

\smallskip
PROOF. Let $X$ and $Y$ be a left and a right normed modules. By the adjoint associativity (cf.,
e.g.,~\cite[Ch.VI.3]{he2}), there is a bijection (actually isometric isomorphism) between ${\bf
h}_{\bf mod-A}(Y,\bb(X,E))={\bf h}_{\bf mod-A}(Y,\bb(?_l,E)(X))$ and ${\bf h}_{\bf
A-mod}(X,\bb(Y,E)={\bf h}_{\bf A-mod}(\bb(?_r,E)(Y), X)$, natural in $X$ and $Y$. The rest is
clear. $\Box$

\medskip
{\bf Example 4.7.} Let us deduce from Proposition 4.5 the description of cofree Banach modules,
contained, sometimes in an equivalent formulation, in~\cite[Ch.III.1]{he1} or~\cite[Ch.VII.1]{he2},
and a similar description of cofree normed modules.

Consider the rigged category $({\bf A-mod},\sq)$, providing `traditional' projective and injective
modules; see Examples 2.4 and 4.1. Take an arbitrary normed space $E$ and consider the diagram
$$
\xymatrix@C+20pt { {{\bf A-mod}^{\bf o}} \ar[r]^{\sq^{\bf o}}\ar[d]_{{\bb(?_l,E)}}
& {\bf Nor}^{\bf o} \ar[d]^{\bb(?_l,E)} \\
{\bf mod-A} \ar[r]^{\sq} & {{\bf Nor} }    },%\eqno(2.2)
$$
\noindent where we use the notation $\sq$ for the obvious forgetful functor in the bottom line, and
the notation ${\bb(?_l,E)}$ for the functor, defined above for $A$-modules and, in particular, for
%an arbitrary $A$ and hence for its specialization in the case of
normed spaces. Take this diagram as
% in the capacity of
(4.1) and the
functor ${\bb(?_l,E)}$ and its specialization for normed spaces as
%in the capacity of
$\Phi$ and $\Delta$, respectively. Because of Proposition 4.6, we see that the conditions of Proposition 4.5
are satisfied.

Evidently, $A$ is a free right normed module with $\co$ as its base normed space. Therefore
Proposition 4.5 implies that the module ${\bb(?_l,E)}(A)$ is a free object in ${\bf A-mod}^{\bf o}$
with respect to the rigged category $({\bf A-mod}^{\bf o},\Box^{\bf o})$, and its base object in ${\bf
Nor}^{\bf o}$ is ${\bb(?_l,E)}(\co)$. Since $\bb(\co, E)=E$, this means that every normed
space $E$ is the cobase space of a  cofree  left normed $A$-module, namely of $\bb(A,E)$ with the
outer multiplication $(a\cd T)(b):=T(ba); a,b\in A, T\in\bb(A,E)$. We see that every
cofree module in ${\bf A-mod}$ has, up to a topological isomorphism, the form $\bb(A,E)$ for a
suitable $E$.

A similar argument works in the `completed' context. Cofree left Banach modules are the
%turn out to be the
same $\bb(A,E)$, only this time $E$ runs the category of Banach spaces.

\medskip
But our main concern here is the cofreedom in the rigged category $({\bf A-mod_1},\circledast)$
from Example 4.3. In what follows, $\circledast$-cofree objects in ${\bf A-mod_1}$ will be called
\emph{metrically cofree normed $A$-modules}. Make the following simple observation.

\medskip
{\bf Proposition 4.8.} \emph{Suppose that all data of Proposition 4.5 are given, and, in addition,
the functor $\Delta$ has a right inverse functor $\nabla:{\cal L}_1\to{\cal L}_2$. Then, if a
rigged category $({\cal K}_2,\sq)$ is freedom-loving, the same is true for $({\cal K}_1,\boxdot)$,
and the free object in ${\cal K}_1$ with the base $M$ is $\Phi\nabla(M)$. }  $\Box$

\medskip
Consider the so-called asterisk functors $(^*_l):{\bf A-mod_1}\to{\bf mod-A_1}$ and $(^*_r):{\bf
mod-A_1}\to{\bf A-mod_1}$, taking a normed module to its dual and a morphism to its adjoint (that
is, particular cases of functors ${\bb(?_l,E)}$ and ${\bb(?_r,E)}$ with $E:=\co$).

\medskip
{\bf Theorem 4.9.} \emph{The rigged category $({\bf A-mod_1},\circledast:{\bf A-mod_1}\to{\bf
Set}^{\bf o})$ is cofreedom-loving, and the cofree module with the cobase set $M$ is $\bb(A,l_\ii(M)$
with the outer multiplication, defined by $[a\cd T](b):=T(ba);a,b\in A, T\in\bb(A,E)$.}

\smallskip
PROOF. Consider the diagram
$$
\xymatrix@C+20pt{  {{\bf A-mod_1}^{\bf o}} \ar[r]^{\:\:\circledast^{\bf o}}\ar[d]_{(^*_l)}
& {\bf Set} \ar[d]^{\id_{\bf Set}} \\
{{\bf mod-A_1}} \ar[r]^{\bigcirc} & {{\bf Set} }  },%\eqno(2.2)
$$
\noindent where the top line is the rig, dual to that in the formulation, and the bottom line is
the obvious `right module version' of the rig in Example 2.7. As the obvious right-module version
of Theorem 2.18, the rigged category $({\bf mod-A_1},\bigcirc)$ is freedom-loving. Moreover, the
right free module with the base set $M$ is $l_1^0(M)\mmp A$ with the outer multiplication, well
defined by $(x\ot b)\cd a:=x\ot ba;a\in A, b\in A, x\in l_1^0(M)$. Taking this into account, we
easily see that the conditions of Proposition 2.27 are satisfied, if we set $\Phi:=(^*_r),
\Psi:=(^*_l)$ and $\Delta:=\Upsilon:=\nabla:=\id_{\bf Set}$. Therefore the rigged category $({\bf
A-mod_1}^{\bf o},\circledast^{\bf o})$ is freedom-loving, and the free object with the base $M$ is
$(^*_r)[l_1^0(M)\mmp A]$. The latter, by the adjoint associativity, is isometrically isomorphic
to the left module $\bb(A,[l_1^0(M)]^*)$, that is to $\bb(A,l_\ii(M))$
%with the outer multiplication, indicated
in the formulation.  $\Box$

%It remains to notice that $({\bf A-mod_1}^{\bf o},\circledast^{\bf o})$ is the rigged category, dual to $({\bf A-mod_1},\circledast)$.

\medskip
A similar description of cofree modules holds, up to obvious modifications, in the case of
completed algebras and modules.

\section{ Metrically projective and free operator spaces}

 Here we pass from `classical' to quantum functional analysis. We shall freely use
 the original, called `matrix' or `coordinate', approach to its concepts and results, presented in
 books~\cite{efr}~\cite{pau}~\cite{pis}~\cite{blem}. The only terminological change is that, speaking about
 what is called in these books `abstract
operator space' and `operator space structure', we shall say `quantum space' and `quantum norm',
respectfully (avoiding the protean adjective `operator'). Thus a \emph{quantum norm} on a linear
space $E$ is a sequence of norms, defined, for every $n\in{\B N}$, on the space $M_n(E)$ of
$n\times n$ matrices with entries in $E$ and satisfying the axioms of Ruan.
%(in their coordinate presentation).
A \emph{quantum space} is a linear space, equipped with a quantum norm.

If $\va:F\to E$ is an operator, its \emph{$n$-amplification} is the operator $\va_n:M_n(E)\to
M_n(F)$, taking the matrix $(a_{kl})$ to $(b_{kl}:=\va(a_{kl}))$. An operator $\va$
between quantum spaces is called \emph{completely bounded}, if $\sup\{\|\va_n\|;n\in{\B N}\}<\ii$;
we call this supremum \emph{completely bounded norm of $\va$} and denote it by $\|\va\|_{cb}$. The
same $\va$ is called \emph{completely contractive, completely isometric,  completely isometric
isomorphism, completely coisometric} or \emph{completely strictly
coisometric}, if the operator $\va_n$  is contractive, respectively, isometric, isometric isomorphism etc. for all $n$.

The space of completely bounded operators between quantum spaces $E$ and $F$, which is a quantum
space in its own right~\cite[p. 45-46]{efr}, is denoted by $\cb(E,F)$.

If $E$ is a quantum space, we identify $M_m(M_n(E))$ with $M_{mn}(E)$ and thus make $M_n(E)$ the
quantum space as well. Note that $\|\va_n\|_{cb}=\|\va\|_{cb}$.

 In what follows,  ${\bf QNor_1}$, respectively, ${\bf QNor}$, denotes the category where objects are
  quantum normed spaces (not necessarily complete), and morphisms are completely contractive operators,
respectively, all completely bounded operators. The `completed' versions of these categories are
denoted by ${\bf QBan_1}$ and ${\bf QBan}$.

Let $A$ be an algebra, endowed with a quantum norm. We call it \emph{quantum algebra}, if the 
bilinear operator of multiplication is completely contractive in the sense of~\cite[p. 126]{efr}. A 
module over a quantum algebra is called \emph{quantum module}, if the bilinear operator of outer 
multiplication is completely contractive. The category of quantum modules and completely 
contractive morphisms, respectively, all completely bounded morphisms will be denoted ${\bf 
QA-mod_1}$, respectively, ${\bf QA-mod}$. Thus ${\bf QNor_1}={\bf Q\co-mod_1}$, and ${\bf 
QNor}={\bf Q\co-mod}$. 

As `quantum' versions of Definitions 1.1 and 1.2, we can define relatively projective and
topologically projective quantum modules. They are projective objects of ${\bf QA-mod}$ with
respect to the rig $\sq:{\bf QA-mod}\to{\bf QNor}$ or $\sq:{\bf QA-mod}\to{\bf Set}$,
where $\sq$ is the appropriate forgetful functor. It is not hard to show that the first rigged
category (contrary to the second) is freedom-loving, and its free object with the base $E$ is $A\ot_{op}E$,
where $\ot_{op}$ is the symbol of the non-completed operator-projective tensor
product of quantum spaces~\cite[p. 124]{efr}. The relevant outer multiplication is defined
similarly to the `classical' case (cf. Section 2).

\bigskip
However, here we are mostly interested in the following kind of projectivity.

\medskip
{\bf Definition 5.1}. A quantum space $P$ is called \emph{metrically projective}, if, for every
completely strictly coisometric operator $\tau:F\to E$ between quantum spaces and an arbitrary
completely bounded operator $\va:P\to E$, there exists a completely bounded operator $\psi$ such
that it is a lifting of $\va$ across $\tau$, and $\|\psi\|_{cb}=\|\va\|_{cb}$.

\smallskip
In a concise form, $P$ is metrically projective, if the morphism functor $\cb(P,?):{\bf
QNor}\to{\bf QNor}$ preserves the property of an operator to be a strict coisometry.

%\medskip
%If, in the given definition, we shall consider completely coisometric operators instead of
%completely strictly coisometric operators we shall come to the spaces, called in~\cite{ble} (just)
%projective. In the present context, it would be natural to call these quantum spaces {\it extremely
%projective.} (Why `extremely', see Section 1) We shall discuss them in the next section.

\medskip
Similarly, with obvious modifications, one can define a notion of a metrically projective
quantum module over a quantum algebra.

In what follows, we restrict ourselves, just for the sake of the simplicity of our presentation,
with the case of quantum spaces ( = quantum $\co$-modules). All subsequent constructions and
results can be extended to the case of general $A$; cf. also the end of the section.

\medskip
Now we suggest a rigged category that enables us to study the metric projectivity by means of the
freedom.
Consider the covariant functor $\bigodot:{\bf QNor_1}\to{\bf Sets}$, acting as follows. It takes a
quantum space $E$ to the set $\textsf{X}_{n=1}^\ii\bigcirc_{M_n(E)}$, the cartesian product of closed
unit balls of the normed spaces $M_n(E)$ (see above). Thus, elements of $\bigodot(E)$ are sequences
$(v_1,\dots v_n,\dots)$ where $v_n\in\bigcirc_{M_n(E)}$. As to morphisms, our functor takes a
completely contractive operator $\va:G\to E$ to the map $\bigodot(\va):\bigodot(G)\to\bigodot(E):
(u_1,\dots,u_n,\dots)\mt(\va_1(u_1),\dots,\va_n(u_n),\dots)$. Evidently, such a functors is a rig.
%faithful. Thus it is a rig.
Note the obvious

\medskip
{\bf Proposition 5.2}. \emph{A completely contractive operator $\tau:G\to E$ between quantum spaces
is an admissible morphism with respect to the rigged category $({\bf QNor_1},\bigodot)$ if, and
only if it is a completely strict coisometry.} $\Box$

\medskip
{\bf Corollary 5.3}. \emph{$\bigodot$-projective quantum spaces are exactly metrically projective
quantum spaces}.

\bigskip
Denote by ${\cal N}_n; n\in{\B N}$ the spaces of trace-class (they say also nuclear) operators on
$\co^n$, endowed with the trace-class norm. We recall that, by virtue of the Schatten/von Neumann
Theorem, there exist the isometric isomorphisms ${\cal N}_n\to(\bb(\co^n))^*$, taking an operator
$b$ to the functional $f:a\mt tr(ba); a\in\bb(\co^n)$. (Here and thereafter $tr(\cd)$ means trace).
This allows us to equip all ${\cal N}_n;n\in{\B N}$ with the quantum norm, induced from the
respective quantum dual spaces; see, e.g.,~\cite[p. 41]{efr}.

Fix, for a time, a quantum space $E$. The existence  of a complete isometric isomorphism between 
the quantum spaces $M_n(E)$ and $\cb({\cal N}_n,E)$, which is the main contents of the next 
proposition, was  indicated by Blecher~\cite[p. 26]{ble}. According to him, it follows from the 
results of~\cite{blep}). 

Recall that, as a linear space, $M_n(E)$, is identified with $M_n\ot E$.
Besides, as a quantum space, $M_n:=M_n(\co)$ is identified with $\bb(\co^n)$. Thus
${\cal N}_n=M_n^*$. Conversely, because of~\cite[Prop. 3.2.1]{efr} and $\dim{\cal N}_n<\ii$,
we have $({\cal N}_n)^*=M_n$.

\medskip
{\bf Proposition 5.4.} {\it There exists a completely isometric isomorphism $\iota_n^E:M_n(E)
\to\cb({\cal N}_n,E)$, well defined by taking an elementary tensor $a\ot x\in M_n\ot E$ to the
operator $b\mt tr(ab)x; a\in M_n=\bb(\co^n), b\in{\cal N}_n, x\in E$. Moreover, $\iota_n^E$ is natural on
$E$; that is, in the detailed form, for every completely bounded operator $\va:G\to E$ between
quantum spaces, there exists the commutative diagram
$$
\xymatrix@C+20pt{ M_n(G) \ar[r]^{\iota_n^G}\ar[d]_{\va_n}
& \cb({\cal N}_n,G) \ar[d]^{\va_n^\bullet} \\
M_n(E) \ar[r]^{\iota_n^E} & \cb({\cal N}_n,E) }\eqno(5.1)
$$
\noindent where $\va_n^\bullet$ is induced by $\va$ (that is, it acts by the rule $\psi\mt\va\psi$).}

\smallskip
PROOF. Denote by ${\cal I}$ the canonical embedding of $E$ into $E^{**}$ and consider the following
chain of quantum spaces and operators:
$$
M_n(E)\stackrel{\iota_1}{\longrightarrow} M_n(E^{**})\stackrel{\iota_2}{\longrightarrow}
\cb(E^*,M_n)\stackrel{\iota_3}{\longrightarrow}\cb({\cal N}_n,E^{**}).
$$
Our operators are as follows:

-- $\iota_1$ is the $n$-amplification of ${\cal I}$, and hence it is a complete isometry together
with the latter~\cite[Prop. 3.2.1]{efr}.

-- $\iota_2$ is the completely isometric isomorphism, participating in the definition of the
quantum dual space~\cite[p. 41]{efr}.

-- $\iota_3$ acts by taking $\va:E^*\to M_n=({\cal N}_n)^*$ (see above) to $\psi:f\mt\beta$, where
$\beta:g\mt[\va(f)]; f\in({\cal N}_n)^*,\beta\in (E^*)^*, g\in E^*$. By the rule of the so-called
quantum adjoint associativity~\cite[p. 128]{efr}, it is a completely isometric isomorphism.

Consequently, the composition of these operators, denoted by $\iota_0:M_n(E)\to\cb({\cal
N}_n,E^{**})$, is a complete isometry. But, taking elementary tensors in $M_n\ot E$, we easily see
that $\iota_0=\iota_4\iota_n^E$, where the operator $\iota_4:\cb({\cal N}_n,E)\to\cb({\cal
N}_n,E^{**})$ is induced by ${\cal I}$. Since both $\iota_0$ and $\iota_4$ are complete isometries
(the latter again by~\cite[Prop. 3.2.1]{efr}) the same is necessarily true for $\iota_n^E$.

Besides, the identification of (${\cal N}_n)^*$ and $M_n$ obviously implies that every operator of
rank 1 within $\cb({\cal N}_n,E)$ acts, for some $a\in M_n=\bb(\co^n)$ and $x\in E$ as $b\mt tr(ab)x$,
and hence coincides with $\iota_n^E(a\ot x)$. Since $\dim{\cal N}_n<\ii$, it follows that
$\iota_n^E$ is a surjection. Hence it is a complete isometric isomorphism.

As to the last claim, concerning the diagram (5.1), we immediately verify it on elementary tensors
in $M_n(G)$. $\Box$

\medskip
In what follows, we rely heavily on the known fact (cf. Blecher~\cite[p. 23]{ble}) that the
category ${\bf QNor_1}$ (as well as ${\bf A-mod_1}$ in Section 2) \emph{admits coproducts.}

Following~\cite{ble}, we denote the coproduct object of the family $E_\nu;\nl$ of quantum spaces by 
$\oplus_1\{E_\nu;\nl\}$. The universal property of this quantum space was explicitly stated in the 
book of Pisier~\cite[p. 52]{pis}, who uses the notation $l_1\{E_\nu;\nl\}$. See also Blecher/Le 
Merdy~\cite[p. 26-27]{blem}. We only note that all these authors write about Banach (= complete) 
operator spaces, whereas we prefer to speak here about the non-complete case. However, the both 
cases are quite similar: the Banach coproduct space is the completion of the coproduct space that 
we use now. 

\medskip
{\bf Remark 5.5.} We shall use only the very fact of the existence of coproducts in ${\bf QNor_1}$,
and not their explicit construction. Nevertheless, for the convenience of the reader, we shall
recall one of possible ways to display them. As a linear space, $\oplus_1\{E_\nu;\nl\}$ is
the algebraic direct sum $\oplus\{E_\nu;\nl\}$; coproduct injections are the natural embeddings of direct summands.

To introduce a quantum norm, that is a `classical' norm in every space $M_n(\oplus\{E_\nu; 
\\ \nl\})$, satisfying the Ruan axioms, we consider the index set $\Upsilon$, consisting of all 
possible pairs $(H_\Upsilon,\{\va_\nu:E_\nu\to\bb(H_\Upsilon);\nl\})$, where $H_\Upsilon$ is a 
Hilbert space (the same for all $\nu$) and $\va_\nu;\nl$ are completely contractive operators. 

Take $u\in M_n(\oplus\{E_\nu;\nl\})$.
Since the latter space coincides, up to a linear isomorphism, with $\oplus\{M_n(E_\nu);\nl\}$, our
$u$ can be treated as the sum $\sum_\nu u_\nu$, where $u_\nu\in M_n(E_\nu)$, with finite number of
non-zero summands. Now consider $\bb(H_\Upsilon)$ with its standard quantum norm and set
$$
\|u\|:=\sup\{\|\sum_\nu(\va_\nu)_n(u_\nu)\|\},
$$
where the supremum is taken over all families $\{\va_\nu;\nl\}$, belonging to pairs in $\Upsilon$.
All needed properties can be easily checked.

%A shorter, but somewhat more sophisticated way is to consider the natural duality between $l_1$-sum
%of $E_\nu$ and $l_\ii$-sum of $E_\nu^*$ (as normed spaces). The latter space has a natural quantum
%norm of the  $l_\ii$-sum of quantum dual spaces. So, we  endow the former space by  the quantum
%norm, generated by the mentioned duality (cf~\cite[Ch. 8.1]{he?}).

\medskip
Combining what was said about coproducts with Proposition 2.12, we obtain

\medskip
{\bf Corollary 5.6.} \emph{The $\oplus_1$-sum of an arbitrary family of metrically projective
quantum spaces is  metrically projective}.

%\medskip
%The statement above, concerning the extreme projectivity (to be precise, its `complete' version),
%was proved in~\cite[Prop. 3.6]{ble}) with the use of the explicit construction of the
%$\oplus_1$-sum.

\medskip
To move further, we need two preparatory statements. The following one belongs to the realm of
linear algebra.

\medskip
{\bf Proposition 5.7.} {\it Suppose that $E,F,G$ are linear spaces, $\dim F=\dim G=n$, and
$e'_1,\dots e'_n, e''_1,\dots e''_n$, where $n\in{\B N}$, are linear bases in $F$ and $G$ respectfully.
Set $u:=\sum_{k=1}^ne''_k\ot e'_k\in G\ot F$. Then, for every $v\in G\ot E$ there exists a unique operator
$\va:F\to E$ such that $[\id_G\ot\va](u)=v$.}

\smallskip
PROOF. We know that $v$ can be presented as $\sum_{k=1}^ne''_k\ot x_k$ with uniquely determined
$x_k\in E$. On the other hand, for every operator $\psi:F\to E$ we have
$\id_G\ot\psi(u)=\sum_{k=1}^ne''_k\ot\psi(e'_k)$. From this, using the linear independence of
$e''_1,\dots e''_n$, we obtain that the desired $\va$ is the only operator, taking $e'_k$ to $x_k;
k=1,\dots,n$. $\Box$

\medskip
Let us concentrate on the special case of Proposition 5.4, when $E$ is ${\cal N}_n$. We have a
completely isometric isomorphism $\iota_n^{{\cal N}_n}:M_n({\cal N}_n)\to\cb({\cal N}_n,{\cal N}_n)$.
Further, we distinguish the special element  ${\bf I}_n$ in $M_n({\cal N}_n)$, uniquely defined by
$i_n^{{\cal N}_n}({\bf I}_n)=\id_{{\cal N}_n}$. (Actually, ${\bf I}_n=\sum_{i,j=1}^ne_{ij}\ot
e_{ji}$, where $e_{ij}$ is the elementary matrix with 1 as `$ij$-th' entry. But we do not need this
observation). We have, of course, $\|{\bf I}_n\|=1$.

\medskip
{\bf Proposition 5.8}. {\it For every quantum space $E$, $n$ and $x\in\bigcirc_{M_n(E)}$ there is a
unique operator $\va:{\cal N}_n\to E$ such that $\va_n:M_n({\cal N}_n)\to M_n(E)$ takes ${\bf I}_n$
to $x$. Moreover, $\va$ is completely contractive.}

\smallskip
PROOF. The first assertion is an obvious particular case of Proposition 5.7, with ${\cal N}_n$ as
$G$ etc. As to the second one, again taking ${\cal N}_n$ as $G$ and looking at (5.1) and at $\va\in
\cb({\cal N}_n,E)$ we see that $\va=\va_n^\bullet(\id_{{\cal N}_n})=i_n^E\va_n({\bf
I}_n)=i_n^E(x)$. Therefore, since $i_n^E$ is isometric, $\|\va\|_{cb}=\|x\|\le1$. $\Box$

\medskip
Finally, we are able to display $\bigodot$-free objects, referred in what follows as
\emph{metrically free quantum spaces}. Denote by $({\cal N}_\ii,\{{\textsf i}^n; n=1,2,\dots\})$
the coproduct of the family $\{{\cal N}_n; n=1,\dots\}$ in ${\bf QNor_1}$. In other words,
$$
{\cal N}_\ii:={\cal N}_1\oplus_1{\cal N}_2\oplus_1\cdots\oplus_1{\cal N}_n\oplus_1\cdots.
$$
Moreover, ${\textsf i}^n:{\cal N}_n\to{\cal N}_\ii$ obviously are coretractions in ${\bf QNor_1}$,
and the same is true with their amplifications $({\textsf i}^n)_n:M_n{\cal N}_n\to M_n({\cal
N}_\ii)$. In particular, these operators are complete isometries.
%Therefore, by Corollary 5.1, the respective amplifications $({\textsf i}^n)_n:\bb_n{\cal
%N}_n\to\bb_n{\cal N}_\ii$ are also coretractions, and in particular, complete isometries.
Set, for every $n$, ${\bf I}^n:=({\textsf i}^n)_n({\bf I}_n)$; we see that $\|{\bf I}^n\|=1$.

\medskip
{\bf Theorem 5.9.} {\it The rigged category $({\bf QNor_1},\bigodot)$ is freedom-loving. Moreover,

(i) the metrically free quantum space with a one-point base, say  $\{t\}$, is ${\cal N}_\ii$.
%The respective universal arrow $\al_{\{t\}}$ takes $t$ to the sequence  $({\bf I}^1,\dots,{\bf
%I}^n,\dots)$ in $\prod\{\bp_{\bb_n{\cal N}_\ii};n=1,2,\dots\}$.

%\medskip
(ii) for a set $M$, the  metrically free quantum space with the base $M$ is the
coproduct space of the family of copies of ${\cal N}_\ii$, indexed by points of
$M$. }

\smallskip
PROOF. (i) Our task is to find, for every quantum space $E$, a bijection
$$
{\cal I}_E:{\bf h}_{\bf Set}(\{t\},\bigodot(E))\to{\bf h}_{\bf Nor_1}({\cal N}_\ii,E),\q\q\eqno(5.2)
$$
\noindent natural in $E$. Take a
%n element from the first of the written sets, that is a
map $\va^0:\{t\}\to\bigodot(E)$. It sends $t$ to some sequence $(x_1,\dots,x_n,\dots);\; x_n\in
M_n(E)$.
By Proposition 5.8, for every $n\in{\B N}$ there exists a unique completely contractive operator
$\va^n:{\cal N}_n\to E$ such that $(\va^n)_n{\bf I}_n=x_n$. By the universal property of coproducts
(cf. Section 2), there exists a unique completely contractive operator $\va:{\cal N}_\ii\to E$ such
that $\va{\textsf i}_n=\va^n$. Thus we obtain a well defined map ${\cal I}_E:\va^0\mt\va$ between our sets of
morphisms.

Further, suppose that a
contractive operator $\va:{\cal N}_\ii\to E$, is given. Assigning to $\va$ the sequence
$(\dots,x_n,\dots)\in\bigodot(E)$, where $x_n:=\va_n({\bf I}^n)$, and then the map from $\{t\}$ to
$\bigodot(E)$, taking $t$ to that sequence, we obtain a map ${\cal J}_E$ from the second of sets in
(5.2) into the first. Since $\va_n({\bf I}^n)=\va_n({\textsf i}^n)_n({\bf I}_n)=(\va^n)_n({\bf
I}_n)$, we easily check that ${\cal I}_E$ and ${\cal J}_E$ are mutually inverse maps.
%Then the element $\va_n({\bf I}^n)
%Then we have $\va_n({\bf I}^n)=\va_n({\textsf i}^n)_n({\bf I}_n)=(\va^n)_n({\bf I}_n)=x_{n}$, and
%this exactly means that $\bigodot(\va)\al_{\{t\}}=\va^0$. The rest is clear.

Thus ${\cal I}_E$ is a bijection. Evidently, it is natural on $E$.

\smallskip
(ii) Immediately follows from (i) and Corollary 2.17. $\Box$

%\medskip
%{\bf Remark 5.10.} In the framework of the so-called non-coordinate approach to the operator space
%theory (see~\cite{he4}), there is a temptation to consider another rig. We mean a functor ${\bf
%Nor_1}\to{\bf Set}$, taking a quantum space $E$ to the closed unit ball of a {\it single} space,
%denoted by ${\cal F}E$ and called the amplification of $E$ within that approach (see {\it idem}, p.
%1). Indeed, such a rig provides projective objects that are, like in the case of $\bigodot$,
%exactly metrically projective quantum spaces. However, this rig behaves much worse than $\bigodot$, when
%we turn to freedom: it provides no free objects save {\bf 0}.

%\medskip
%Now we came to the point, where rigged categories $({\bf QNor_1},\bigodot)$ and $({\bf QNor_1},
%\%slash)$ behave essentially differently. Namely, the argument, similar to that in the proof of
%Pr%position 1.?, shows that the second rigged category (like $({\bf QNor_1},\bp)$ and $({\bf
%QNor_1}, \circledcirc)$), has no free objects save 0. We conjecture that there is no rig of ${\bf
%QNor_1}$ that would provide  extremely projective quantum spaces and at the same time have non-zero free spaces.

\medskip
Let us call a quantum space \emph{nuclear-composed} (or trace-class-composed), if it is of the form
$\oplus_1\{E_\nu;\nl\}$, where each of summands is ${\cal N}_n$ for some $n\in{\B N}$.
Incidentally, Theorem 5.9 shows that such a space is metrically free if, and only if the
cardinality of summands ${\cal N}_n$ with fixed $n$ is the same for all $n$.

%Theorem 5.1 immediately provides the respective ``metric'' counterpoint: combining this theorem
%with Propositions 1.? and 1.?, we obtain
Combining Theorem 5.9 with Proposition 2.11, we obtain

\medskip
{\bf Corollary 5.10.} {\it (i) Every quantum space is completely strict coisometric image of a
nuclear-composed quantum space.

(ii) A quantum space is metrically projective if, and only if it is a retract in ${\bf QNor_1}$ of
a nuclear-composed quantum space (equivalently, it is a direct summand of such a space with the
respective natural projection of completely bounded norm 1). In particular, all ${\cal N}_n;
n\in{\B N}$ are metrically projective. }

\medskip
This corollary is a `metric counterpart' of the results of Blecher~\cite[Prop.3.1 and Thm
3.9]{ble} on the extremely projective quantum spaces (see Section 6 below). The argument in~\cite{ble} is
based on work with quantum duals of mentioned coproducts.

% We recall that Blecher~\cite[Prop.3.1 and Thm
%3.9]{ble} (working with completed spaces) has shown that every quantum space is a quotient ( =
%completely coisometric image) of an ${\cal N}$-built quantum space, and that a given quantum space
%is extremely projective if, and only if it is a so-called almost direct summand of an ${\cal
%N}$-built quantum space. In equivalent form (cf. `near-retracts' in~\cite[]{he9}) this means that,
%for every $\e>0$, there is a completely contractive operator from the latter onto the former space,
%possessing a right inverse $\rho$ with $\|\rho\|_{cb}<1+\e$.

Of course, the metric projectivity of ``bricks'' ${\cal N}_n$ immediately follows from the 
identification of $M_n(E)$ and $\cb({\cal N}_n, E)$, provided by Proposition 5.4 (cf.~\cite[Prop. 
3.7]{ble}). 

%\medskip
%As to the last assertion, there is another way to prove it. Namely, the metric, as well as the
%extreme projectivity of  ${\cal N}_n$
% can be immediately derived from the identification  $\bb_nE$ and $\cb({\cal N}_n, E)$,
%provided by Proposition 5.6 (cf.~\cite[Prop. 3.7]{ble}). Besides, combining  p.(ii) of the previous
%corollary with the above-mentioned Blecher's characterization of extremely projective quantum
%spaces, we see that every metrically projective quantum space is extremely projective. We do not
%know whether the converse is true.

\medskip
Theorem 5.9 can be easily generalized from quantum spaces to quantum modules
% over quantum algebras
(see the beginning of the section). The relevant rig is acting from ${\bf QA-mod_1}$ to ${\bf Set}$ and,
similarly to $\bigodot$ above, takes an $A$-module $X$ to $\textsf{X}_{n=1}^\ii\bigcirc_{M_n(X)}$
The resulting rigged category is again freedom-loving, and the respective free modules are the
$\oplus_1$-sums of families of copies of the quantum $A$-module $A\ot_{op}{\cal N}_\ii$.
% (cf. the beginning of this section).

In another direction, one can study injective and cofree quantum spaces (and modules) with respect to the quantum
version of the rig $\circledast$ from Section 4. In particular, cofree quantum spaces will happen to be
quantum $l_\ii$-sums of quantum spaces $\bb(\co^n); n=1,2,\dots$.

\section{ Asymptotic structure in categories and extreme projectivity}

Definitions and results, concerning the extreme projectivity,
necessarily have, so to say, `asymptotic nature'. Here we suggest the general scheme,
including such a kind of projectivity. Throughout the section ${\cal K}$ is an arbitrary category.

\medskip
{\bf Definition 6.1.} Let $\{{\B J}_\nu; \nl\}$ be a family of natural transformations of the
identity functor on ${\cal K}$ into itself. (Recall that in the detailed form this means that for
every object $X$ in ${\cal K}$ a morphism ${\B J}_\nu^X:X\to X$ is given, and for every $\nl$ and a
morphism $\va:X\to Y$ the diagram
$$
\xymatrix@C+20pt{ X \ar[r]^{\va}\ar[d]_{{\B J}_\nu^X}
& Y \ar[d]^{{\B J}_\nu^Y} \\
X \ar[r]^{\va} & Y }
$$
\noindent is commutative). Such a family is called \emph{asymptotic structure on} ${\cal K}$, if
it satisfies the following two conditions:

(i) for every $\nl$ and $X$, the morphism ${\B J}_\nu^X$ is a bimorphism.

(ii) for every $\nl$, there are $\lm,\mu\in\Lambda$ such that ${\B J}_\nu={\B J}_\lm{\B J}_\mu$.

\medskip
{\bf Definition 6.2.} A triple, consisting of a category, a rig of the latter and an asymptotic
structure on this category, is called \emph{ asymptotic category}.

\medskip
From now on, we suppose that we are given an asymptotic category  \\
$({\cal K},\sq:{\cal K}\to{\cal L},\{{\B J}_\nu; \nl\})$.

\medskip
{\bf Definition 6.3.} A morphism $\va:X\to Y$ in ${\cal K}$ is called \emph{permitted},
(with respect to the given asymptotic structure)
%$\{{\B J}_\nu; \nl\}$,
if there exists $\nl$ and a morphism $\widetilde\va:X\to Y$ such that $\va=\widetilde\va{\B
J}_\nu^X$ (or, equivalently, $\va={\B J}_\nu^Y\widetilde\va$).

\medskip
{\bf Definition 6.4.}
A morphism $\tau:Y\to X$ in ${\cal K}$ is called \emph{asymptotically admissible epimorphism}, if,
for every $\nl$, there exists a morphism $\rho_\nu:\sq(X)\to\sq(Y)$ in ${\cal L}$ such that we have
$\sq(\tau)\rho_\nu\sq({\B J}_\nu^X)=\sq({\B J}_\nu^X)$. The family $\{\rho_\nu;\nl\}$ will be
referred as {\it asymptotically right inverse} to $\sq(\tau)$.

\medskip
The word `epimorphism' above is justified. Indeed, if we have $\va\tau=\psi\tau$, then \\
$\sq(\va)\sq(\tau)\rho_\nu\sq({\B J}_\nu^X)=\sq(\psi)\sq(\tau)\rho_\nu\sq({\B J}_\nu^X)$.
Hence $\sq(\va)\sq({\B J}_\nu^X)=\sq(\psi)\sq({\B J}_\nu^X$). Since ${\B J}_\nu^X$ is epi (Definition 6.1), and
% together with the faithfulness of 
$\sq$ is faithful, we have $\va=\psi$.

%Note an obvious

\medskip
{\bf Proposition 6.5.} \emph{Every admissible epimorphism is asymptotically admissible. As a 
corollary, if our asymptotic category, being considered just as a rigged category, is 
freedom-loving, then every object in ${\cal K}$ is a range of an asymptotically admissible morphism 
with a free domain.} 

\smallskip
PROOF. We just take one-point $\Lambda$ and the `genuine' right inverse to $\sq(\tau)$ as the
unique morphism $\rho_\nu$. The rest is clear. $\Box$

\medskip
Finally, we suggest

\medskip
{\bf Definition 6.6.} An object $P$ in ${\cal K}$ is called \emph{asymptotically projective} if,
for every asymptotically admissible epimorphism $\tau$, every permitted morphism $\va:P\to X$ has
a lifting in ${\cal K}$ across $\tau$.

\medskip
{\bf Example 6.7.} (`classical'). Consider the rigged category in Example 2.7. Set $\Lambda:=(0,1)$ 
and, for every $t\in(0,1)$ and a normed $A$-module $X$, set ${\B J}_t^X:X\to X:x\mt tx$. This is an 
asymptotic structure on ${\bf A-mod_1}$. Indeed, the properties of a natural transformation, 
required in Definition 6.1, as well as (i), are obvious, and (ii) is valid because ${\B J}_t^X=[{\B 
J}_{\sqrt{t}}^X]^2$. We obtain an asymptotic category with ${\bf A-mod_1}$ as ${\cal K}$ and $\bp$ 
as $\sq$. Permitted morphisms are clearly those with norm $<1$. Asymptotically admissible 
epimorphisms are exactly coisometries. Indeed, if we are given a coisometry $\tau:Y\to X$, we 
choose as $\rho_t:\bp_X\to\bp_Y$ the map, taking a vector $x$ to $y$, where $y$ is: 

-- an arbitrary vector in $\bp_Y$ with $\tau(y)=x$, provided $\|x\|\le t$

-- an arbitrary vector in $\bp_Y$ (without any condition) otherwise.

It is easy to verify that asymptotically projective modules turn out to be extremely projective in the
sense of Definition 1.3.

%Finally, asymptotically projective modules are those $P$ with the following property: for every
%coisometric $A$-module morphism $\tau:Y\to X$, every bounded $A$-module  morphism $\va:P\to X$ and
%every $\e>0$, there exists a  bounded $A$-module  morphism $\psi:P\to Y$, such that (i) $\psi$ is a
%lifting of $\va$ across $\tau$, and (ii) $\|\psi\|<\|\va\|+\e$. With an accordance to~\cite{??}, we
%call modules with the indicated property \emph{extremely projective}. (The term is chosen because
%the morphism functor ${\bf h}_A(P,?):{\bf A-mod}\to{\bf Nor}$ preserves the property of a morphism
%to be, in general-categorical sense, an extreme epimorphism; cf~\cite[]{clm}).

\medskip
{\bf Example 6.8} (`quantum'). Take the rigged category $({\bf QNor_1},\bigodot)$ from Section 5. 
Again, set $\Lambda:=(0,1)$ and, for $t\in(0,1)$ and a quantum space $E$, set ${\B J}_t^E:E\to 
E:x\mt tx$. Like in the previous example, we obtain an asymptotic category, this time with ${\bf 
QNor}$ as ${\cal K}$ and $\bigodot$ as $\sq$. Now an operator $\va$ is a permitted morphism exactly 
when $\|\va\|_{cb}<1$, and asymptotically admissible epimorphisms are complete coisometries. 
Indeed, if $\tau:G\to E$ is a complete coisometry, we choose as $\rho_t:\textsf{X}\{M_n(E); 
n=1,2,\dots\}\to\textsf{X}\{M_n(G); n=1,2,\dots\}$ the map, taking a sequence $(\dots,x_n,\dots)$ 
to $(\dots,y_n\dots)$, where $y_n$ is: 

-- an arbitrary vector in $\bp_{M_n(G)}$ with $\tau(y_n)=x_n$, provided $\|x_n\|\le t$

-- an arbitrary vector in $\bp_{M_n(G)}$ (without any condition) otherwise.)

It is easy to verify that a quantum space $P$ is asymptotically projective if, and only if it has 
the following property: for every completely coisometric operator $\tau:G\to E$, every completely 
bounded operator $\va:P\to E$ and every $\e>0$, there exists a completely bounded operator 
$\psi:P\to G$, such that (i) $\psi$ is a lifting of $\va$ across $\tau$, and (ii) 
$\|\psi\|_{cb}<\|\va\|_{cb}+\e$. These quantum spaces were introduced in~\cite{ble} and called 
there (just) projective; we shall call them 
%, for the same reason as in the previous example, 
\emph{extremely projective}.

\medskip
This asymptotic structure, as well as Corollary 6.13(ii) below, can be easily extended from quantum spaces to quantum modules over an arbitrary quantum algebra.

Return to the general scheme. Speaking about free objects of an asymptotic category and free-loving asymptotic categories,
 we mean
 %, of course, a free object with respect to 
 the underlying rigged category.

\medskip
{\bf Proposition 6.9.} \emph{Every free object of an asymptotic category}
\emph{is asymptotically projective}.

\smallskip
PROOF. Let ${\bf Fr}(M)$ be a $\sq$-free object with the base object $M$. Suppose we are given an 
asymptotically admissible epimorphism $\tau:Y\to X$ with respective approximately right inverse 
$\rho_\nu;\nl$ and a permitted morphism $\va:P\to E$; thus the latter 
%as we remember, 
has the form ${\B J}_\mu^X\widetilde\va$ for some $\mu\in\Lambda$ and $\widetilde\va:P\to X$.

Remembering about bijections, participating in Definition 2.10, set   \\ 
$\al:={\cal I}_{\bf Fr}(M)^{-1}(\id_{\bf Fr(M)}):M\to\sq({\bf Fr}(M))$, $\psi^0:=\rho_\mu\sq(\va)\al:M\to\sq(Y)$ 
and finally $\psi:={\cal I}_Y(\psi^0):{\bf Fr}(M)\to Y$. Since the mentioned bijections are natural in the second 
argument, we have
$$
\tau\psi=\tau[{\cal I}_Y(\psi^0)]={\cal I}_X[\sq(\tau)\psi^0]={\cal I}_X[\sq(\tau)\rho_\mu\sq(\va)\al]=
{\cal I}_X[\sq(\tau)\rho_\mu\sq({\B J}_\mu^X\widetilde\va)\al]=
$$
$$
{\cal I}_X[\sq({\B J}_\mu^X)\sq(\widetilde\va)\al]=
{\cal I}_X[\sq(\va)\al]=
\va {\cal I}_{{\bf Fr}(M)}[\al]=\va{\cal I}_{{\bf Fr}(M)}{\cal I}_{{\bf Fr}(M)}^{-1}[\id_{{\bf Fr}(M)}]=
\va. \:\Box
$$

%\smallskip
%PROOF. Let $F$ be a $\sq$-free object with the base object $M$and a universal arrow$\al:M\to\sq(F)$,
%$\tau:Y\to X$ an asymptotically admissible epimorphism, $\rho_\nu;\nl$ the
%respective approximately right inverses to $\tau$, $\va:P\to E$ a permitted morphism; the latter,
%as we remember, has the form ${\B J}_\mu^X\widetilde\va$ for some $\mu\in\Lambda$ and
%$\widetilde\va:P\to X$.
%
%Set $\psi^0:=\rho_\mu\sq(\va)\al:M\to\sq(Y)$. Then the freedom of $F$ provides a unique morphism
%$\psi$ in ${\cal K}$ such that $\sq(\psi)\al=\psi^0$. We have
%$$
%\sq(\tau\psi)\al=\sq(\tau)\rho_\mu[\sq({\B J}_\mu^X)\sq(\widetilde\va)]\al=
%\sq({\B J}_\mu^X)\sq(\widetilde\va)]\al=
%\sq(\va)\al.
%$$
%It follows that $\tau\psi=\va$. $\Box$

\medskip
{\bf Proposition 6.10.} \emph{Let $P$ be an asymptotically projective object, $\tau:Y\to X$ an
asymptotically admissible epimorphism, and $\va:P\to X$ a permitted morphism. Then the set of all liftings of
$\va$ across $\tau$ contains a permitted morphism}.

\smallskip
PROOF. The assumption on $\va$, combined with Definition 6.1(ii), provides $\lm,\mu\in\Lambda$
and $\widetilde\va:P\to X$ such that $\va={\B J}_\lm^X{\B J}_\mu^X\widetilde\va$. Since ${\B
J}_\mu^X\widetilde\va$ is permitted, the assumption on $P$ provides $\widetilde\psi:P\to Y$ such
that $\tau\widetilde\psi={\B J}_\mu^X\widetilde\va$. Set $\psi:=\widetilde\psi{\B J}_\lm^Y$. We
have
$$
\tau\psi=\tau\widetilde\psi{\B J}_\lm^Y={\B J}_\mu^X\widetilde\va{\B J}_\lm^Y={\B J}_\lm^X{\B J}_\mu^X\widetilde\va=\va.
\; \Box
$$

\medskip
{\bf Definition 6.11.} A morphism $\s:U\to V$ in ${\cal K}$ is called \emph{asymptotic
retraction}, if, for every $\nl$, there exists a morphism $\zeta_\nu:V\to U$ such that
$\s\zeta_\nu={\B J}_\nu^V$. An object $V$ in ${\cal K}$ is  called \emph{asymptotic retract} of
an object $U$, if there exists an asymptotic retraction from $U$ into $V$.

\medskip
Consider the asymptotic category of Example 6.7. It is easy to see that a contractive morphism
of normed $A$-modules
% $\s:U\to V$ 
is an asymptotic retraction if, and only if for every $\e>0$ it has a right inverse module morphism with norm $<1+\e$. 
Thus an asymptotic retraction is exactly a near-retraction in the terminology of~\cite{he7}.

In a similar way, one can easily describe asymptotic retractions in the context of the asymptotic
category of Example 6.8. Namely, a completely contractive operator $\s:U\to V$, where $U,V$
are now quantum spaces, is an asymptotic retraction if, and only if for every $\e>0$ it has a
right inverse completely bounded operator with completely bounded norm $<1+\e$. We shall call again such
an operator a near-retraction (this time in ${\bf QNor}$). One can easily see that $V$ is an asymptotic
retract of $U$ if, and only if it almost direct summand of $U$ in the terminology of~\cite{ble}.

\medskip
{\bf Proposition 6.12.} \emph{Suppose that our asymptotic category is
%, being considered just as a rigged category, is 
freedom-loving. Then an object $P$ in ${\cal K}$ is asymptotically projective if, and only if it is 
an asymptotic retract of a free object}. 

\smallskip
PROOF. `Only if' part. By Propositions 2.11(i) and 6.5, there exist a free object $F$ and an 
asymptotically admissible epimorphism $\tau:F\to P$.
Take $\nl$. Since ${\B J}_\lm^P={\B J}_\lm^P\id_P:P\to P$ is a permitted morphism, there is a
%exists its
lifting, say $\zeta_\nu$, across $\tau$. The rest is clear.

`If' part. Let $\s:F\to P$ be an  asymptotic retraction with a free domain, $\tau:Y\to X$ an 
asymptotically admissible morphism, and $\va:P\to X$ a permitted morphism. Then it easily follows 
from the diagram in Definition 6.1 that $\va\s:F\to X$ is permitted as well. Hence, by Propositions 
6.9 and 6.10, there exists $\chi:F\to Y$ such that $\tau\chi=\va\s$ and $\chi={\B 
J}_\nu^Y\widetilde\chi$ for some 
 $\widetilde\chi:F\to Y$ and $\nl$. Set $\psi:=\widetilde\chi\zeta_\nu$, where $\zeta_\nu$ is such that
 $\s\zeta_\nu={\B J}_\nu^P$ (see Definition 6.11). We have
$$
{\B J}_\nu^X\tau\psi={\B J}_\nu^X\tau\widetilde\chi\zeta_\nu=\tau{\B J}_\lm^Y\widetilde\chi\zeta_\nu=
\tau\chi\zeta_\nu=\va\s\zeta_\nu=\va{\B J}_\nu^P={\B J}_\nu^X\va.
$$
But the morphism ${\B J}_\lm^X$ is mono (see Definition 6.1(i)); hence $\tau\psi=\va$. $\Box$

\medskip
Combining this with what was said about the asymptotic categories in Examples 6.7 and 6.8, and also
with Theorems 2.8 and 5.9 together with Proposition 6.12, we immediately have

\medskip
{\bf Corollary 6.13.} \emph{ (i) A normed module $P$ over a normed algebra $A$ is extremely
projective if, and only if it is a near-retract of a module of the form $A\ot_p l_1^0(M)$, where
$M$ is a set. The same is true for a Banach module over a Banach algebra after replacing `\;$\ot_p$'
by `\;$\prtpet$' and $l_1^0$ by $l_1$.}

\emph{ (ii)} (Blecher~\cite[Thm. 3.10]{ble}) \emph{A quantum space $P$ is extremely projective if, 
and only if it is a near-retract of a coproduct in ${\bf QNor_1}$ of 
%quantum spaces of trace-class operators 
some family of spaces  ${\cal N}_n; n\in{\B N}$. 
%of finite-dimensional trace-class spaces. 
The same is true after replacing 
the words `quantum space' by `complete quantum space' and 
 ${\bf QNor_1}$ by ${\bf QBan_1}$.}

\ed
\begin{thebibliography}{999}

\bibitem{ble}
D.~P.~Blecher. The standard dual of an operator space, Pacific J. of Math. \\ v. 153, No. 1 (1992) 15-30.
\bibitem{blem}
D.~P.~Blecher, C.~Le~Merdy. {\it Operator algebras and their modules}.  Clarendon Press. Oxford.
2004.
\bibitem{blep}
D.~P.~Blecher, V.~I.~Paulsen. Tensor products of operator spaces, J. Funct. Anal., 99, No. 2 (1991)
262-292 .
\bibitem{clm}
J.~Cigler, V.~Losert, P.~Michor. {\it Banach modules and functors on categories of Banach spaces}.
Marcel Dekker, New York, 1979.
\bibitem{efr}
E.~G.~Effros, Z.-J.~Ruan. {\it Operator spaces}.  Clarendon Press. Oxford. 2000.
\bibitem{faith}
C.~Faith. {\it Algebra: Rings, modules, and categories I.} Springer-Verlag. Berlin. 1981. 
\bibitem{gro}
A.~Grothendieck. Une caracterisation vectorielle-metrique des espaces $L^1$, Canadian J. Math., 7
(1955) 552-561.
\bibitem{groe}
N.~Groenbaek. Lifting problems for normed spaces. Preprint.
\bibitem{gus}
E.~V.~Gusarov. General-categorical frame-work for topologically free normed modules. Prepront. 
2012. 
\bibitem{he0}
A.~Ya.~Helemskii. On the homological dimensions of normed modules over Banach algebras, Mat. Sb.
(N.S.) 81 (123) (1970) 430-444; transl.: Math. USSR Sb. 10 (1970) 399-411.
\bibitem{he1}
A.~Ya.~Helemskii. {\it The Homology of Banach and Topological Algebras}. Kluwer, Dordrecht, 1989.
\bibitem{he2}
A.~Ya.~Helemskii. {\it Banach and locally convex algebras}. Clarendon Press. Oxford. 1993.
\bibitem{he3}
A.~Ya.~Helemskii. {\it Lectures and exercises on functional analysis}. AMS, Providence, R.I., 2005.
\bibitem{he4}
A.~Ya.~Helemskii. {\it Quantum functional analysis. Non-Coordinate Approach}. AMS, Providence,
R.I., 2010.
\bibitem{he5}
A.~Ya.~Helemskii. Extreme flatness of normed modules and Arveson-Wittstock type theorems, J.
Operator Theory, 64:1 (2010) 101-112.
\bibitem{he6}
A.~Ya.~Helemskii. Metric version of flatness and Hahn-Banach type theorems
 for normed modules over sequence algebras, Studia Math. 206 (2), 2011.
\bibitem{he7}
A.~Ya.~Helemskii. Extreme version of projectivity for normed modules over sequence algebras. 
arXiv:1104.2463 v 1 [math FA]. 
%\bibitem{fpo}
%M.~Frank, V.~I.~Paulsen. Injective and projective Hilbert $C^*$-modules, and $C^*$-algebras of
%compact operators, arXiv:math/0611349v2 [math.OA] 18 Feb 2008.
\bibitem{heb}
G.~K\"othe. Hebbare lokalkonvexe R\"aume, Math. Annalen 165 (1966) 181-195.
%\bibitem{loe}
%H.~L\"owig. \"Uber die Dimension linearer R\"aume. Studia Math. 5 (1934), 18-23.
\bibitem{mc1}
S.~Mac Lane. {\it Homology}.  Springer-Verlag, Berlin, 1967.
\bibitem{mc2}
S. Mac Lane. {\it Categories for the working mathematician}. Springer-Verlag, Berlin, 1971.
\bibitem{pau}
V.~I.~Paulsen. {\it Completely bounded maps and operator algebras}. Cam. Univ. Press, Cambridge, 2002.
\bibitem{pie}
A.~Pietsch. {\it History of Banach spaces and linear operators.} Birkh\"auser. Basel. 2007.
\bibitem{pis}
J.~Pisier. {\it Introduction to operator space theory}. Cam. Univ. Press. Cambridge, 2003.
\bibitem{run}
V.~Runde. {\it Lectures on Amenability}. Springer-Verlag, Berlin, 2002.
\bibitem{sem0}
Z.~Semadeni. {\it Projectivity, injectivity and duality}. Rozprawy Matematyczne, PAN, Warzsawa, 35 
1963. 
\bibitem{sem}
Z.~Semadeni. {\it Banach spaces of continuous functions}. PWN. Warsaw. 1971.
\bibitem{tsal}
Tsalenko, Shul'geifer.


\bibitem{wit}
G.~Wittstock. Injectivity of the module tensor product of semi-Ruan modules, J. Operator Theory, 65:1 (2011) 87-113.

\end{thebibliography}
